\documentclass[preprint,12pt,sort&compress]{elsarticle}

\usepackage{amssymb}

\usepackage{lineno}

\usepackage[total={6.5in, 9in}]{geometry}
\usepackage{amsfonts, amsmath,amsthm}
\usepackage{xcolor}
\usepackage[hyperfootnotes=false]{hyperref}
\hypersetup{
    colorlinks,
    linkcolor={red!50!black},
    citecolor={blue!50!black},
    urlcolor={blue!80!black}
}
\usepackage{graphicx}
\usepackage{booktabs}
\usepackage{mathtools}
\usepackage[scientific-notation=true]{siunitx}
\usepackage{algorithm}
\usepackage{algorithmic}
\usepackage{cleveref}
\floatname{algorithm}{Algorithm} 

\newcommand{\ve}{\varepsilon}

\newcommand{\RR}{\mathbb{R}}

\newcommand{\uscat}{u^{\textrm{scat}}}
\newcommand{\uinc}{u^{\textrm{inc}}}
\newcommand{\utot}{u^{\textrm{tot}}}
\newcommand{\umeas}{u^{\textrm{meas}}}
\newcommand{\cF}{\mathcal{F}}
\newcommand{\pa}{\partial}

\newcommand{\parentheses}[1]{\left(#1\right)}

\newcommand{\norm}[1]{\left\|#1\right\|}

\newcommand{\mz}[1]{\textcolor{red}{\textsf{[MZ: #1]}}}

\newcommand{\revisea}[1]{#1}
\newcommand{\reviseb}[1]{#1}
\DeclareMathOperator*{\argmin}{arg\,min}

\newcommand{\bc}{\boldsymbol{c}}
\newcommand{\bcnn}{\boldsymbol{c}^{\textrm{NN}}}
\newcommand{\bcmeas}{\boldsymbol{c}^{\textrm{meas}}}
\newcommand{\bclsm}{\boldsymbol{c}^{\textrm{LSM}}}
\newcommand{\Ntrain}{N_{\textrm{train}}}
\newcommand{\Nepo}{N_{\textrm{epoch}}}

\newtheorem{remark}{\sffamily Remark}

\begin{document}

\begin{frontmatter}



\title{A Neural Network Warm-Start Approach for the Inverse Acoustic Obstacle Scattering Problem}


\author[inst1]{Mo Zhou}
\ead{mo.zhou366@duke.edu}

\affiliation[inst1]{organization={Department of Mathematics},
            addressline={Duke University}, 
            city={Durham},
            postcode={27708}, 
            state={NC},
            country={USA}}

\author[inst2]{Jiequn Han\corref{cor1}}
\ead{jiequnhan@gmail.com}
\author[inst2]{Manas Rachh}
\ead{mrachh@flatironinstitute.org}
\author[inst3]{Carlos Borges}
\ead{carlos.cardosoborges@ucf.edu}
\cortext[cor1]{Corresponding author.}
\affiliation[inst2]{organization={Center for Computational Mathematics},
            addressline={Flatiron Institute}, 
            city={New York},
            postcode={10010}, 
            state={NY},
            country={USA}}

\affiliation[inst3]{organization={Department of Mathematics},
            addressline={University of Central Florida}, 
            city={Orlando},
            postcode={32816}, 
            state={FL},
            country={USA}}

\begin{abstract}
In this paper, we consider the inverse acoustic obstacle problem for sound-soft star-shaped obstacles in two dimensions wherein the boundary of the obstacle is determined from measurements of the scattered field at a collection of receivers outside the object. One of the standard approaches for solving this problem is to reformulate it as an optimization problem: finding the boundary of the domain that minimizes the $L^2$ distance between computed values of the scattered field and the given measurement data. The optimization problem is computationally challenging since the local set of convexity  shrinks with increasing frequency and results in an increasing number of local minima in the vicinity of the true solution. In many practical experimental settings, low frequency measurements are unavailable due to limitations of the experimental setup or the sensors used for measurement. Thus, obtaining a good initial guess for the optimization problem plays a vital role in this environment. 

We present a neural network warm-start approach for solving the inverse scattering problem, where an initial guess for the optimization problem is obtained using a trained neural network.  We demonstrate the effectiveness of our method with several numerical examples. For high frequency problems, this approach outperforms traditional iterative methods such as Gauss-Newton initialized without any prior (i.e., initialized using a unit circle), or initialized using the solution of a direct method such as the linear sampling method. The algorithm remains robust to noise in the scattered field measurements and also converges to the true solution for limited aperture data. However, the number of training samples required to train the neural network scales exponentially in frequency and the complexity of the obstacles considered. We conclude with a discussion of this phenomenon and potential directions for future research.
\end{abstract}




\end{frontmatter}


\section{Introduction}
The inverse acoustic scattering problem arises naturally, inter alia, in sonar, radar, medical imaging, 
and seismic detection \cite{kuchment2014radon, collins1995nondestructive, engl2012inverse, Ustinov2014, cheney2009fundamentals, beilina2015globally, thanh2015imaging}. 
In this problem, one or several incident waves are used to recover different properties of the domain, such as its shape, its density, and variations of sound velocity within the medium. 

In this work, we consider the inverse acoustic scattering problem for sound-soft obstacles in two dimensions. 
The forward problem for the scattered field $\uscat$ in the time-harmonic setting is given by
\begin{equation} \label{eq:dir_problem}
\begin{cases}
\Delta \uscat+k^2 \uscat = 0, \quad \text{in} ~\mathbb{R}^{2} \setminus \overline{D}, \\
\uscat = - \uinc \, \quad \text{on} ~\partial D \, , \\
\lim_{r\rightarrow \infty} r^{1/2}\left(\frac{\partial \uscat}{\partial r} - ik\uscat\right) = 0 \,, 
\end{cases} 
\end{equation}
where $k$ is the wave number, $\partial D$ is the boundary of the obstacle $D$, and $\uinc$ is the incident field. Without loss of generality, we will consider $\uinc (x) = e^{ik x \cdot d}$, a plane wave with wavenumber $k$ and incident direction $d \in S^1$. Let $x_{j}$, $j=1,2,\ldots N_{t}$, denote a collection of receivers located sufficiently far away from the obstacle, i.e., $\min_{j} \text{dist}(x_{j}, \overline{D}) \gg 2 \pi/k$. We define the forward scattering operator $\cF_{k,d} : \mathcal{B} \to \mathbb{C}^{N_{t}}$, as the scattered field evaluated at the receivers, i.e., 
\begin{equation}
\label{eq:fscat}
     \cF_{k,d}(\pa D) = \umeas_{k,d} \, ,
\end{equation} 
where the $j$th component of $\umeas_{k,d}$ is given by $\uscat_{k,d}(x_{j})$. 
Here $\mathcal{B}$ denotes the set of non-intersecting simply connected curves.
Given measurements at the receivers from one or multiple incident waves, the inverse problem corresponding to the forward scattering problem in~\cref{eq:fscat} is to recover the shape of the obstacle $\pa D$. In particular, let $d_{\ell}$, $\ell=1,2,\ldots N_{d}$, denote the incident directions, and let $\umeas_{k,d_{\ell}}$ denote the corresponding measurements at the receivers. Then, the inverse obstacle scattering problem seeks the shape that minimizes the following objection function 
\begin{equation}
\label{eq:optim}
\widetilde{\pa D} = \argmin_{\pa D} \sum_{\ell=1}^{N_{d}} \| \umeas_{k,d_{\ell}} - \cF_{k,d_{\ell}}(\pa D)  \|^2 \, .
\end{equation}

The inverse obstacle scattering problem is inherently non-linear and ill-posed without additional constraints. The ill-posedness can in part be attributed to Heisenberg's uncertainty principle for waves, which states that one cannot stably recover more resolution than half the wavelength of the incoming wave \cite{chen1997inverse}. Thus, in order to obtain a sharp reconstruction of the obstacle, i.e., resolve the high curvature regions or small features of the obstacle, the inverse problem would have to be solved at high frequencies. 
On the other hand, with increasing $k$, the optimization problem becomes increasingly non-convex, with an increasing number of local minima near the global minimum, and the local set of convexity in the vicinity of the global minimum is shrinking as $O(1/k)$~\cite{borges2020inverse}.

In many practical settings, the combination of these difficulties associated with the inverse problem is addressed via multifrequency measurements --- the low frequency data 
enables the solution at each frequency to be close enough to the global minimum, 
while the high frequency measurements allow for the stable reconstruction of small features 
of the obstacle~\cite{chen1995recursive,chen1997inverse,bao2005inverse,bao2012shape,sini2012inverse,chaillat2012faims,borges2015inverse,bao2015inverse}. 
However, scattered data measurements are typically unavailable for low frequencies due to 
limitations of the experimental setup or the type of sensors being used. 
In this regime, solving the single frequency problem at the lowest frequency 
suffers from all the difficulties of solving the inverse problem in~\cref{eq:optim}.

One of the standard approaches for solving the inverse problem is to treat it as a non-linear optimization problem as stated in~\cref{eq:optim}, and use iterative methods such as steepest descent or Gauss-Newton methods. Examples of these iterative methods include \cite{hanke1995landweber, hettlich1999second, hohage1997logarithmic, kirsch1993domain, kress2003newton, kress1994quasi, kress1997integral}. A major drawback of these methods is that the initial guess needs to be close enough to the solution to guarantee convergence. As noted above, the task of finding a good initial guess can be particularly challenging at higher frequencies due to the shrinking local set of convexity near the true solution.

On the other hand, direct methods like sampling based methods are also commonly used 
for the solution of the inverse problem. In these approaches, the boundary of the 
object is typically represented as the level set of an indicator function that is 
computed via the solution of a regularized linear system. In certain regimes, 
it has been shown that the indicator function remains bounded only for all interior points 
of the domain as the regularization parameter is varied. Examples of sampling methods include 
the linear sampling method \cite{colton1996simple}, the generalized linear sampling method 
\cite{audibert2014generalized}, the factorization method \cite{kirsch1998characterization}, singular source methods \cite{potthast2000stability, potthast2001point}, probe methods \cite{ikehata1998reconstruction, ikehata1999reconstruction}, and many others \cite{colton1998inverse}. A major advantage of sampling methods is that they do not require any a priori information about the obstacle such as its topology or boundary conditions. However, sampling methods tend to require a large number of measurements of the scattered field from multiple incident directions to obtain reasonable approximations of the shape.


With the development of deep learning (DL) methods over the last decade, there has been an increasing interest in solving the inverse scattering problem with DL. Examples include the Switchnet structure developed by Khoo and Ying~\cite{khoo2019switchnet}, the BCR neural network by Fan and Ying \cite{yuwei2022solving}, the use of radial basis function neural networks~with Gaussian kernel by Rekanos in \cite{rekanos2002neural}, the back-propagation scheme of Wei and Chen in \cite{wei2018deep}, and many others~\cite{adler2017solving,chen2019learning,guo2019supervised,sanghvi2019embedding}. For a comprehensive review of this class of methods, we refer the reader to~\cite{chen2020review}. From the knowledge of the authors, the complexity of geometries considered in this work as measured in wavelengths of the incident waves is much larger than most of the examples in other DL approaches considered so far, and thus poses more challenges as explained in the next section.

In \cite{bao2007inverse}, the authors presented a robust inverse obstacle scattering solver as a hybrid of both direct and iterative methods. While their approach was designed for the multifrequency inverse problem, its restriction to the single frequency case is essentially a direct imaging method to provide an initial guess for a Gauss-Newton iterative method. Motivated by this success, we present a DL warm-start approach to initialize a traditional iterative solver. In particular, we construct a neural network to approximate the map from the measurements of the scattered field to the obstacle boundary. The output of the trained model is subsequently refined using a Gauss-Newton approach for the optimization problem~\cref{eq:optim}. 
Our approach falls into the ``Learning-Assisted Objective-Function Approach'' category based on the classification in~\cite{chen2020review}. 

We compare the performance of our method with a Gauss-Newton based iterative approach initialized using a unit disk, a modification of the linear sampling method adapted to far-field measurements of the scattered field (as opposed to far-field pattern measurements), and also with a hybrid approach where the Gauss-Newton based iterative solver is initialized using the linear sampling method.

We restrict our attention to star-shaped obstacles in two dimensions, i.e., the boundary of the obstacle is represented by $\Gamma(t): [0,2\pi) \to \mathbb{R}^{2}$, with 
\begin{equation}
\label{eq:star-shape}
\Gamma(t) = r(t) (\cos{(t)}, \sin{(t)}),
\end{equation}
where $r(t)$ referred to as the radius function, is a positive periodic function with period $2\pi$.
The complicated nature of the set of non-intersecting curves poses a considerable challenge for any inverse obstacle scattering solver. The collection of star-shaped curves, while restrictive, forms a rich family of curves whose complexity can be controlled via the number of Fourier modes representing the radius function. 
This allows us to compare different approaches for the solution of the inverse obstacle scattering problem independent of the difficulties that arise from the representation of generic non-intersecting curves, and allows us to understand the impact of both the frequency of the scattering problem and the Fourier content of the curve on the performance of the various methods.

\revisea{\begin{remark}
A more appropriate definition of the complexity of the curve would be a measure of the near bandlimit of its curvature for an arclength parameterization of the curve. Restricted to star-shaped domains, the maximum Fourier content of the curvature of the curve can be increased in one of two ways: either by increasing the bandlimit of the Fourier series of the radius function where each coefficient has a fixed variance, or by holding the number of terms in the Fourier expansion fixed, and increasing the variance of the individual coefficients. In this work, we opt to work in the former regime.
\end{remark}}

The rest of the paper is organized as follows. In Section~\ref{sec:theoretical_background}, we reformulate the inverse obstacle scattering problem in~\cref{eq:optim} for star-shaped obstacles, and review properties of the optimization landscape as a function of frequency. In Section~\ref{sec:existing_methods}, we review a Gauss-Newton iterative method and an adaptation of the linear sampling method for its solution. In Section~\ref{sec:numerical_method}, we describe our deep learning framework. In Section~\ref{sec:numerical_results}, we demonstrate the efficiency of our approach through several numerical examples and compare our method to four other methods for solving the inverse obstacle scattering problem. 
We conclude with a discussion of these results and directions for future work in Section~\ref{sec:conclusion}.

\begin{remark}
In the interest of brevity, whenever we want to refer to a quantity for all $N_{d}$ incident directions, we will omit the subscript corresponding to the direction. For example, the forward scattering operator for all incident directions is given by $\cF_{k} := [\cF_{k,d_{1}} \, ; \cF_{k,d_{2}} \,; \ldots \cF_{k,d_{N_{d}}}]$,  and the measured data is given by $\umeas_{k} := [\umeas_{k,d_{1}} \, ; \umeas_{k,d_{2}} \,; \ldots \umeas_{k,d_{N_{d}}}] \in \mathbb{C}^{N_{t} N_{d}}$.
\end{remark}

\section{Inverse scattering problem}\label{sec:theoretical_background}
In this section, we recap the inverse obstacle scattering problem in the context of recovery of star-shaped obstacles (specified in \cref{eq:star-shape}) and review relevant properties of the optimization landscape of the inverse problem. 

The inverse obstacle scattering problem is inherently ill-posed due to Heisenberg's uncertainty principle for waves, which states that it is difficult to recover small features of the obstacle with respect to the wavelength of the incident wave. 
A mechanism by which this behavior manifests itself is that there is a large family of boundaries with significantly different curvatures, which result in nearly the same far-field measurements at the receivers, thereby making it difficult to distinguish these boundaries in the inverse problem setup.

\begin{figure}[h!]
    \centering
    \includegraphics[width=0.4\textwidth]{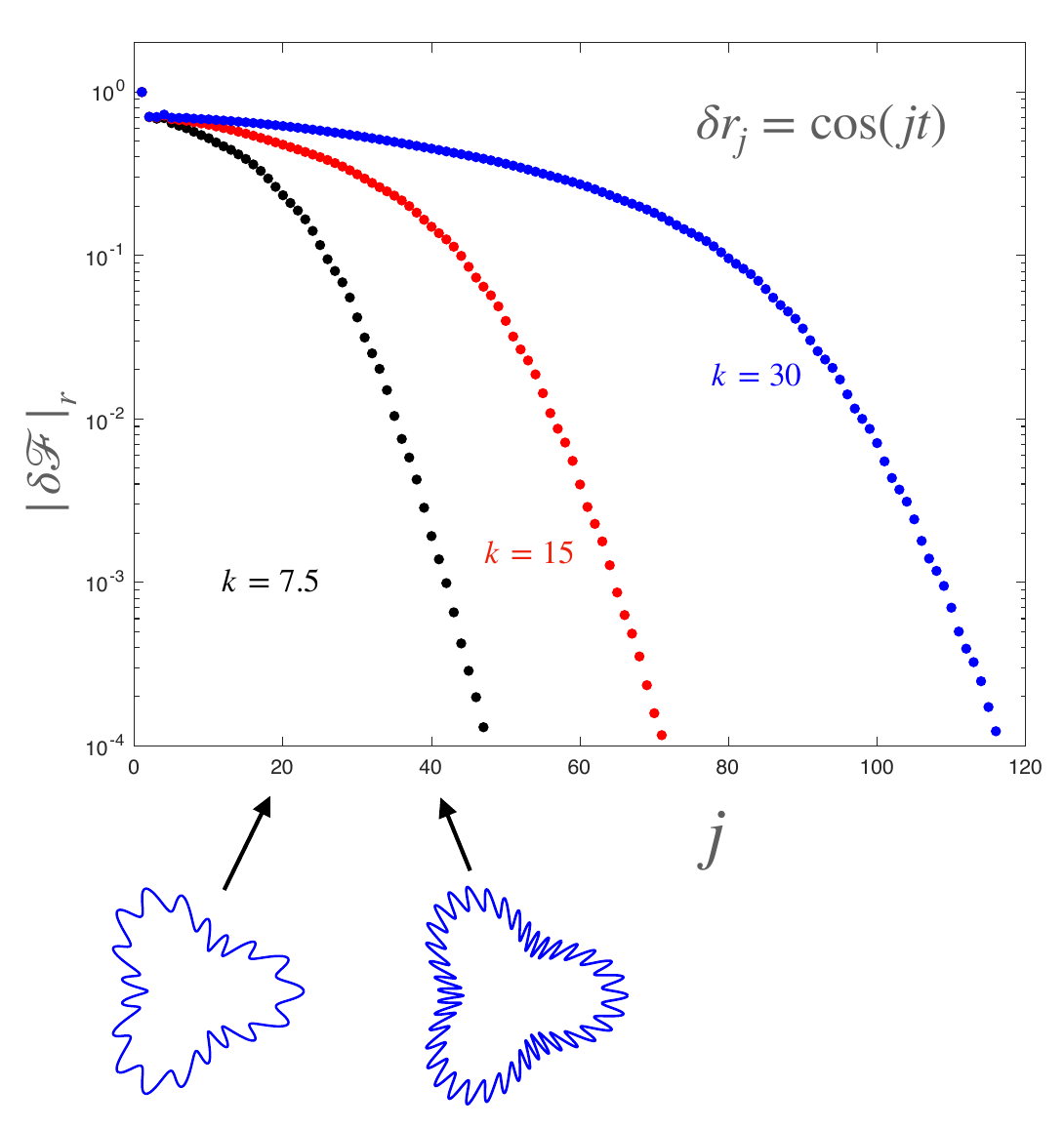}
    \caption{Magnitude of Fr\'echet derivative, $|\delta \mathcal{F}|_{r}$ for $\delta r_{j} := \cos{(j t)}$, scaled by its maximum value over all $j$.} 
    \label{fig:ill-posedness} 
\end{figure}

To illustrate this behavior, consider a star-shaped boundary with $r(t) = 1 + 0.3 \cos{(3t)}$. Suppose that measurements are made at $N_{t} = 200$ equally spaced receivers on the circle of radius $10$ for $N_{d}=200$ incident directions equally spaced in angle. Let $\partial \cF_{k} \delta r$ denote the Fr\'echet derivative of the forward scattering operator $\cF_{k}(\pa D)$ in the direction $\delta r$.  We postpone the details of the computation of the Fr\'echet derivative to Section~\ref{sec:existing_methods}, but briefly, the Fr\'echet derivative measures the sensitivity of the objective function to the boundary perturbation $\delta r$, and each component of the  Fr\'echet derivative $\partial \cF_{k,d_{\ell}} \delta r$ can be obtained via the solution of a related Helmholtz boundary value problem. In Figure~\ref{fig:ill-posedness}, we plot the magnitude of the Fr\'echet derivative for $\delta r_{j} := \cos{(j t)}$ (scaled by the maximum of the norm of the Fr\'echet derivative over all $j$) as a function of $j$ for three different frequencies $k=7.5, 15$, and $30$.

The figure illustrates two properties of the inverse problem. First, for a fixed $k$, the Fr\'echet derivative is $O(1)$ as long as the number of oscillations in the perturbation is less than $O(k)$. As the number of oscillations in the perturbation increases even further, the norm of the Fr\'echet derivative decreases exponentially. 
Thus the measured data for 
$r(t) + \alpha \delta r(t)$ for $j \gg O(k)$ and small enough $\alpha$ would be nearly indistinguishable from 
the measured data for $r(t)$. To remedy this ill-posedness, we constraint the optimization problem so that $r(t)$ is given
by
\begin{equation}
\label{eq:rdef}
 r(t; \bc)=c_{0} + \sum_{m=1}^M \left(c_{m}\cos(mt)+c_{m+M} \sin(mt)\right), 
\end{equation}
with $M = O(k)$, $c_{j}$, $j=0,1,\ldots 2M$, such that $r(t)>0$ for all $t \in [0,2\pi)$, and $\bc = [c_{0}; c_{1} \ldots c_{2M}]$. To summarize, the optimization function takes the form 
\begin{equation}
\label{eq:optim_star}
\tilde{\bc} = \argmin_{\substack{\bc \\ r(t; \bc) > 0 \, \, \forall t\in [0,2\pi)}} \| \umeas_{k} - \cF_{k}(\bc)  \|^2 \, ,
\end{equation}
where in a slight abuse of notation, $\cF_{k}(\bc)$ is the forward scattering operator corresponding to the star-shaped domain with radius function given by $r(t;\bc)$.
\reviseb{\begin{remark}
\Cref{eq:optim_star} avoids ill-posedness by restricting $M = O(k)$, and thus, does not require any further regularization. The condition number of the Fréchet derivative remains bounded throughout the optimization iterations, as illustrated by an example geometry in  Figure~\ref{fig:ill-posedness}.
\end{remark}}

\begin{figure*}[h!]
    \centering
    \includegraphics[width=0.4\textwidth]{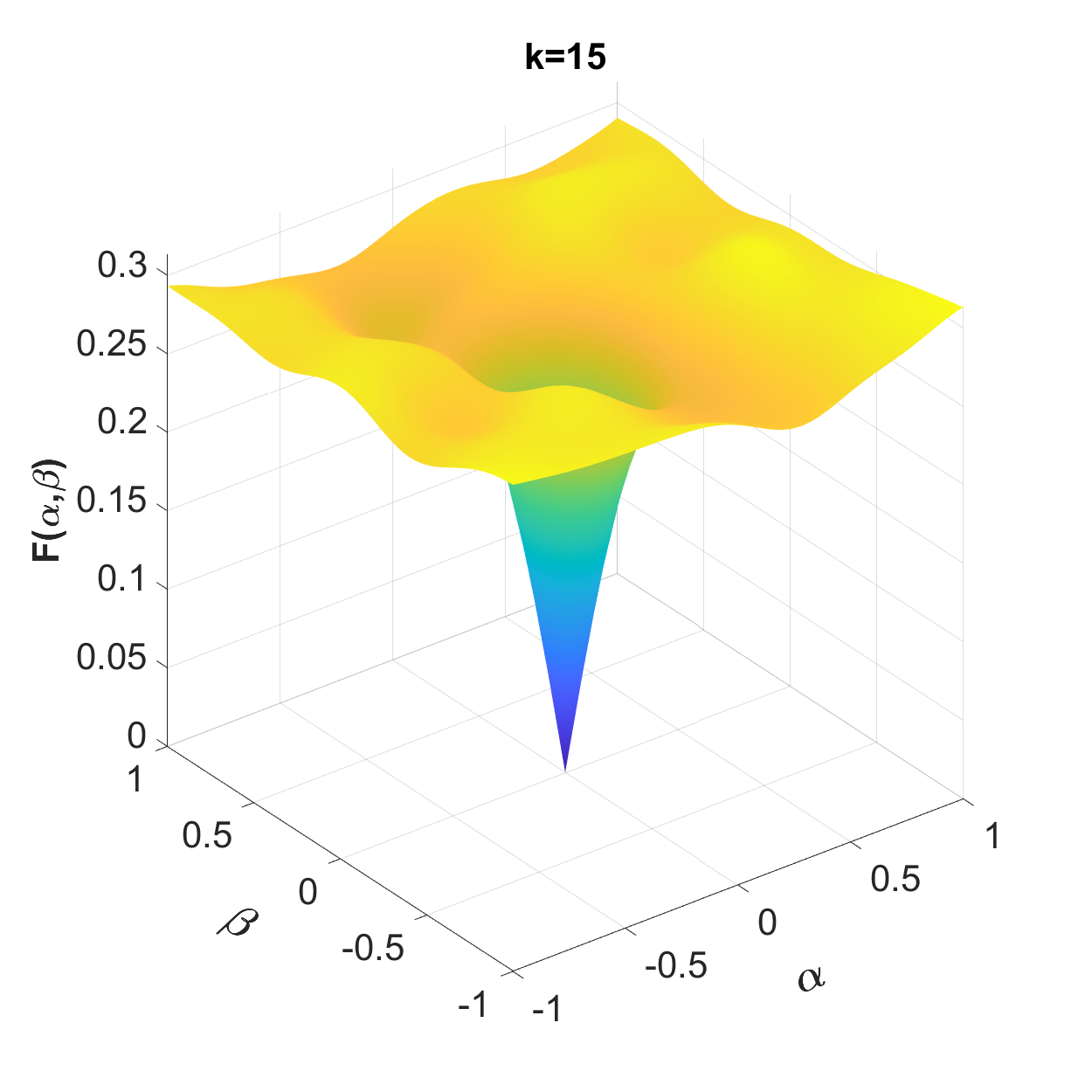}
    \includegraphics[width=0.4\textwidth]{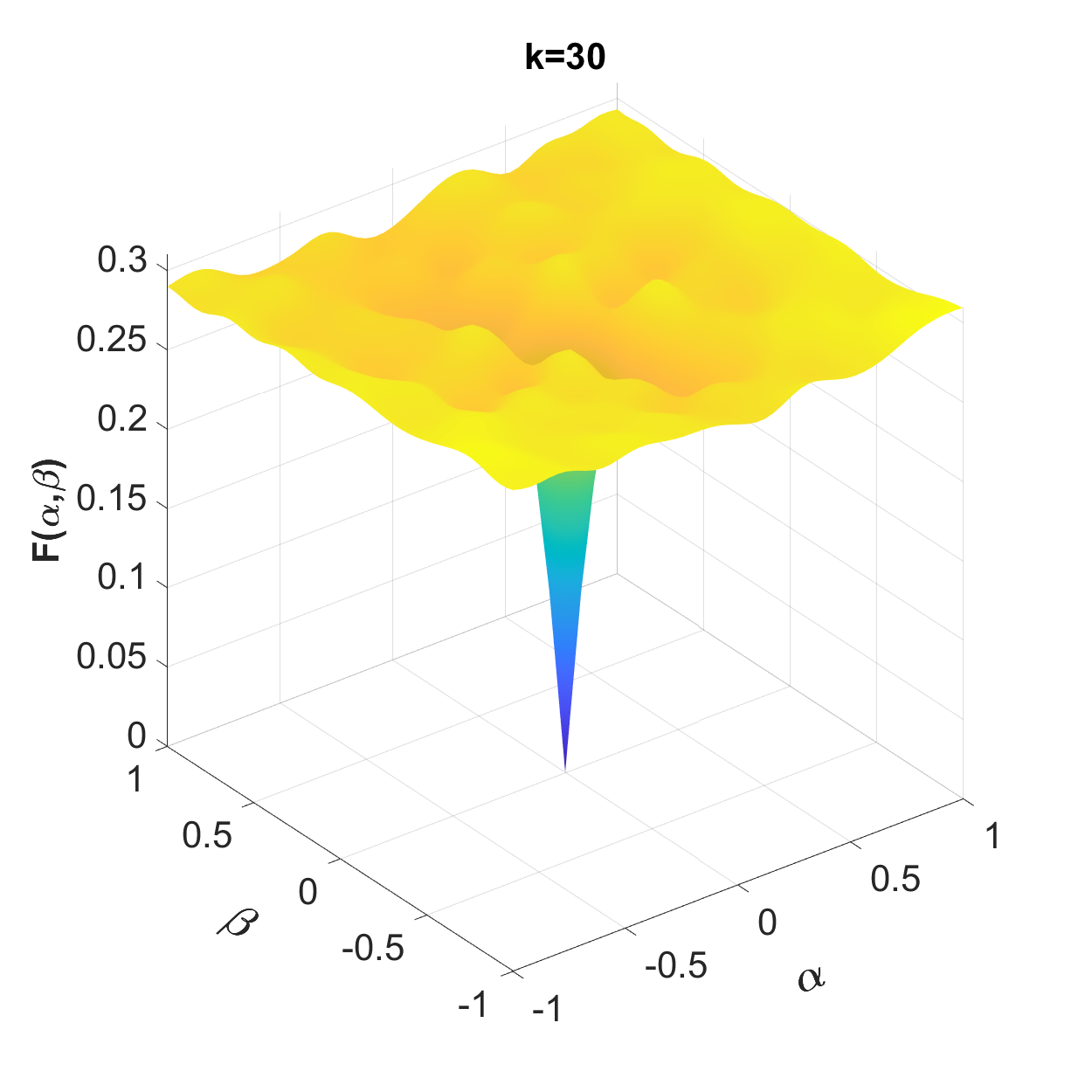}\\
    \includegraphics[width=0.4\textwidth]{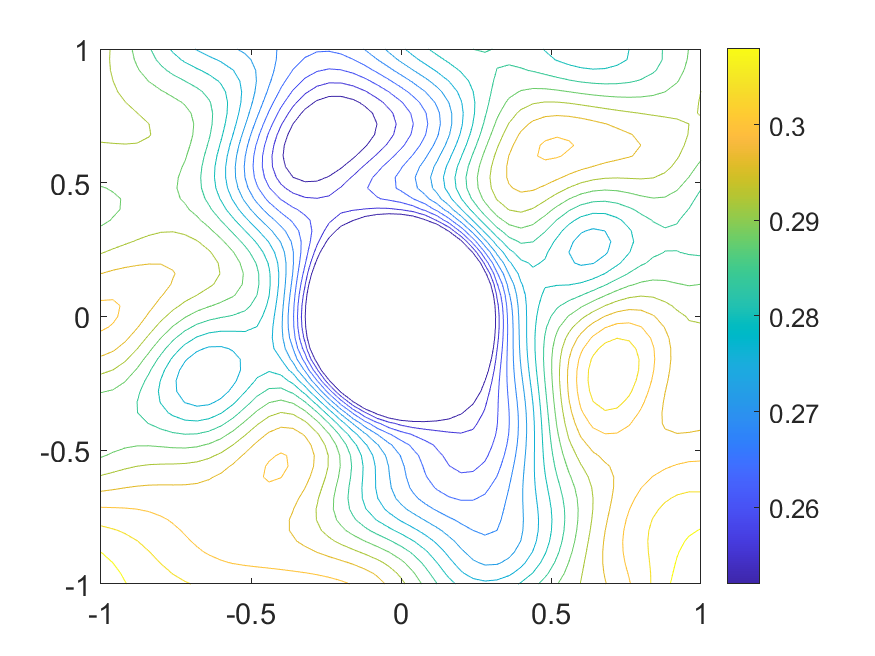}
    \includegraphics[width=0.4\textwidth]{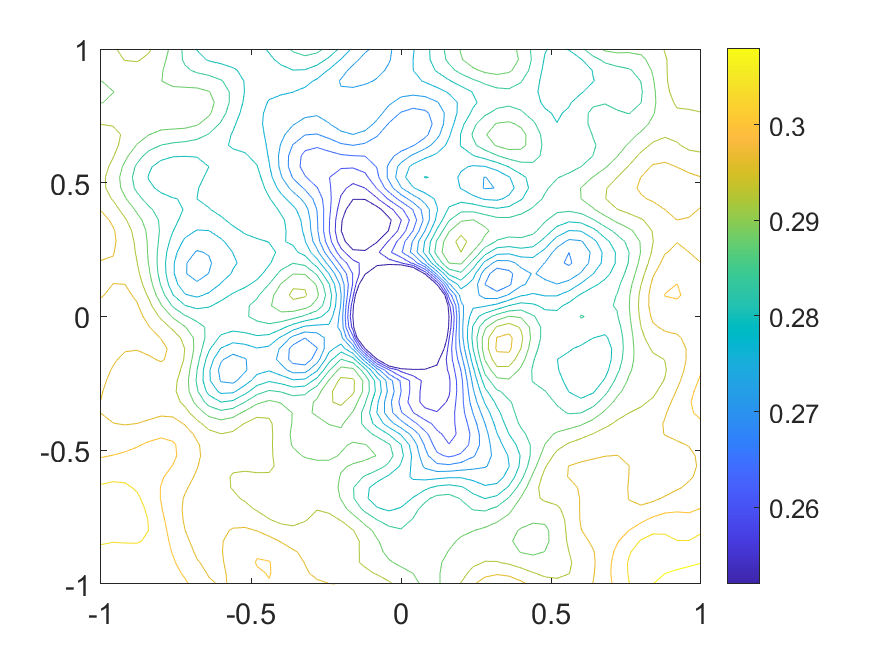}
    \caption{\revisea{Surface plot (top) and contour plot (bottom) of $F(\alpha,\beta)$ defined in~\eqref{eq:local_minima} for $k=15$ (left) and $k=30$ (right).} } 
    \label{fig:landscape} 
\end{figure*}

The inverse obstacle scattering problem is also computationally challenging due to its extreme non-convexity and nonlinearity, and the presence of several local minima in the vicinity of the global minimum.
To illustrate this, we plot the following two-dimensional cross-section of the function being minimized, 
\begin{equation}
\label{eq:local_minima}
F(\alpha,\beta) = \| \cF_{k}(\bc^{\textrm{meas}}) - \cF_{k}(\bc^{\textrm{meas}} + \alpha \bc_{1} + \beta \bc_{2}) \| \, ,
\end{equation}
 for $k=15$, and $k=30$ in Figure~\ref{fig:landscape}. Here $\bc^{\textrm{meas}}$ is such that $r(t;\bc^{\textrm{meas}}) = 1 + 0.3 \cos{(3t)}$, and $\bc_{1}$, and $\bc_{2}$ are random coefficient vectors with $M=5$.

As is evidenced from the figure, the number of local minima increases in the vicinity of the true minimum $(\alpha,\beta) = 0$, as $k$ is increased. It also illustrates that the local set of convexity in the vicinity of the global minimum shrinks like $O(1/k)$. Thus, any iterative method would require an increasingly accurate initial guess for the solution to converge to the global minimum as the frequency is increased.

\section{An iterative and a direct method for the solution of the inverse problem\label{sec:existing_methods}}
We now turn our attention to reviewing two existing approaches for the solution of~\cref{eq:optim_star}. 
In particular, we review a Gauss-Newton iterative approach in~\Cref{sec:Newton} and an adaptation of the linear sampling approach in~\Cref{sec:LSM}.
\subsection{Gauss-Newton method} \label{sec:Newton}
Given a frequency $k$, let $M = O(k)$, and let $\bc^{(0)} \in \mathbb{R}^{2M+1}$ denote an initial guess for a star-shaped domain. In the absence of the positivity constraint, the Gauss-Newton method for the optimization problem~\cref{eq:optim_star}, iteratively updates the solution via
\begin{equation}
\bc^{(j+1)}  = \bc^{(j)} + \delta \bc^{(j)} \, ,
\end{equation}
where $\delta \bc^{(j)}$ is the solution of the least-square problem
\begin{equation}
\label{eq:lsq}
[\, \text{Re}(J) \, ; \text{Im}(J) \,] \cdot \delta \bc^{(j)} = [ \, \text{Re}(\umeas_{k} - \cF(\bc^{(j)})) \, ;  \, \text{Im} (\umeas_{k} - \cF(\bc^{(j)})) \,] \, ,
\end{equation}
with $J \in \mathbb{C}^{(N_{t} N_{d}) \times (2 M+1)}$ being the matrix of Fr\'echet derivatives. \reviseb{Specifically, let $e_{\ell}$ be the $\ell$th coordinate vector in
$\mathbb{R}^{2M+1}$, then the $\ell$th column of $J$ is 
the Fr\'echet derivative in the direction $\delta \bc = e_{\ell}$ and is given by
\begin{equation}
[ \partial \cF_{k, d_{1}} (\bc^{(j)}) e_{\ell}\, ; \,  \partial \cF_{k, d_{2}} (\bc^{(j)}) e_{\ell} \ldots \partial \cF_{k, d_{N_{d}}} ( \bc^{(j)}) e_{\ell}] \, ,
\end{equation} 
where $\partial \cF_{k,d_{m}} (\bc^{(j)})$ is the Fr\'echet derivative with respect to the boundary for the Helmholtz Dirichlet problem with data $\uinc = e^{ik  x\cdot d_{m}}$, which can be computed via the solution of a related Helmholtz boundary value problem, as detailed below.}

To simplify notation, suppose the domain $D$ is star-shaped whose boundary $\pa D$ is parametrized by $\bc$, i.e., $r(t; \bc) = c_{0}  + \sum_{m=1}^{M} (c_{m} \cos{(mt)} + c_{m+M} \sin{(mt)})$. Let $\nu$ denote the outward normal to the boundary $\pa D$. Let $\uscat$ denote the solution to the boundary value problem~\cref{eq:dir_problem}, $\utot = \uscat + \uinc$, and $\cF_{k,d}(r(t; \bc))$ denote the corresponding forward scattering operator.
Let $\delta \bc \in \mathbb{R}^{2M+1}$ parameterize a perturbation to the boundary, i.e., \reviseb{$\delta r(t) = \left(\delta c_0 + \sum_{m=1}^{M} (\delta c_{m} \cos{(mt)} + \delta c_{m+M} \sin{(mt)})\right)\left(\cos(t),\sin(t)\right)$}. Then the Fr\'echet derivative 
$\left(\partial \cF_{k,d} (\bc) \delta \bc \right)_{j} = v(x_{j})$, $j=1,2\ldots N_{t}$, where $v$ is the solution to
\begin{equation}
\begin{aligned}
\Delta v+k^2 v = 0, \quad \text{in} ~\mathbb{R}^{2} \setminus \overline{D}, \\
v = - (\nu \cdot \delta r) \frac{\partial \utot}{\partial \nu}  \, \quad \text{on} ~\partial D \, , \\
\lim_{r\rightarrow \infty} r^{1/2}\left(\frac{\partial v}{\partial r} - ikv\right) = 0 \,. 
\end{aligned}
\end{equation}
\reviseb{A detailed proof for this result can be found in~\cite{potthast1994frechet}. For proof containing the far field pattern, we refer the reader to~\cite{colton1998inverse}.}

Returning back to the Gauss-Newton iteration, note that we split the real and imaginary parts of both the Fr\'echet derivative and the residual vector $\umeas_{k} - \cF(\bc^{(j)})$ in~\cref{eq:lsq} since the update to the curve is parametrized as real numbers, i.e., $\delta \bc^{(j)}\in \mathbb{R}^{2M+1}$, while the measurements of the scattered field are complex. 
Furthermore, recall that the coefficients $\bc^{(j)}\in \mathbb{R}^{2M+1}$ must satisfy the constraint of the corresponding radius function being positive in order to disallow self-intersecting boundaries. If the updated coefficient vector $\bc^{(j+1)}$ does not satisfy the constraint, then we filter the update $\delta \bc^{(j)}$ using a Gaussian filter as
\begin{equation}
\begin{aligned}
\delta c^{(j)}_{m} \to \delta c^{(j)}_{m} \exp{\left( - \frac{m^2}{\sigma^2 M^2} \right)} \, , \quad 0\leq m\leq M \, ,\\
\delta c^{(j)}_{m+M} \to \delta c^{(j)}_{m+M} \exp{\left( - \frac{m^2}{\sigma^2 M^2} \right)} \, , \quad 1\leq m\leq M \, .
\end{aligned}
\end{equation}

We repeat the filtering process up to $10$ times, with $\sigma = 1/10^{(\ell-1)}$, where $\ell$ is the iteration number for the filtering step. If the curve continues to be self-intersecting after $10$ filtering attempts, then we terminate the Gauss-Newton iteration. 
The Gauss-Newton iteration is also terminated if the number of iterations exceeds a maximum iteration count $N_{\textrm{max}}$, or the size of the update $\| \delta \bc^{(j)} \|$ as measured in $\ell^{2}(\mathbb{R}^{2M+1})$ is less than a given update tolerance $\varepsilon_{s}$, or the value of the objective function 
$\| \umeas_{k} - \cF_{k}(\bc^{(j)}) \|$ as measured in $\ell^{2}(\mathbb{C}^{N_{t} \cdot N_{d}})$ is less than a given residual tolerance $\varepsilon_{r}$. The algorithm is summarized in~\Cref{alg}.

\begin{algorithm}%
\caption{Gauss-Newton method for star-shaped obstacles}
\begin{algorithmic}
\REQUIRE Scattered field measurements $\umeas$, an initial guess $\bc^{(0)} \in \RR^{2M+1}$, maximum number of iteration $N_{max}$, update tolerance $\varepsilon_s$, and residual tolerance $\varepsilon_r$.
\STATE{Set $j=0$, $\|\delta\bc^{(j)}\|=2 \epsilon_s$.}
\WHILE{$j<N_{max}$ {\bf and}  $\|{\bf u}^\emph{meas}_k- \mathcal{F}_k(\bc^{j)})\|> \epsilon_r$ {\bf and} $\|\delta\bc^{(j)}\|> \epsilon_s$}
\STATE{Calculate $\cF_k(\bc^{(j)})$ and the Fr\'echet derivative matrix $J$.}
\STATE{Solve $[\text{Re}(J)\,; \text{Im}(J)] \delta \bc^{(j)} = [\text{Re}(\umeas_{k} - \cF_{k} (\bc^{(j)}))\,; \text{Im}(\umeas_{k} - \cF_{k} (\bc^{(j)}))]$.}
\STATE{$\sigma=1$}
\STATE{$\bc^{(j+1)}\leftarrow\bc^{(j)}+ \delta\bc^{(j)}$}
\WHILE{$\min_{t\in[0,2\pi)} r(t,\bc^{(j+1)}) \leq 0$}
\STATE{$\widetilde{\delta c}^{(j)}_{m} \leftarrow \delta c^{(j)}_{m} \exp{\left( - \frac{m^2}{\sigma^2 M^2} \right)} \, , \quad 0\leq m\leq M \, ,$}
\STATE{$\widetilde{\delta c}^{(j)}_{m+M} \leftarrow \delta c^{(j)}_{m+M} \exp{\left( - \frac{m^2}{\sigma^2 M^2} \right)} \, , \quad 1\leq m\leq M \, ,$}
\STATE{$\sigma \leftarrow$ $\sigma/10$}
\STATE{$\bc^{(j+1)}\leftarrow \bc^{(j)}+\widetilde{\delta\bc}^{(j)}$}	
\ENDWHILE
\STATE{$j\leftarrow j+1$}
\ENDWHILE
\end{algorithmic}
\label{alg}
\end{algorithm}

The evaluation of the residual $\umeas - \cF_{k}(\bc^{(j)})$, and the Fr\'echet derivative $J$ at each Gauss-Newton iterate requires solutions to the Helmholtz Dirichlet problem with a fixed boundary but with different boundary data. In this work, we use a standard combined field integral equation representation for the numerical evaluation of both of these quantities, see~\cite{kress1989linear} for example. The discretized linear systems corresponding to the integral equations, and the least square problem in~\cref{eq:lsq} are solved using dense linear algebra methods. For higher frequency problems, the solution to the integral equation could be obtained using fast direct solvers~\cite{chandrasekaran2006fast, chandrasekaran2006fast1, gillman2012direct, greengard2009fast, ho2012fast, martinsson2005fast,bebendorf2005hierarchical, bormhierarchical, borm2003introduction}, while the least square problem could be solved using an iterative approach like LSQR \cite{paige1982lsqr}, for example.

\subsection{The linear sampling method}\label{sec:LSM}
The linear sampling method (LSM) was first introduced by Colton and Kirsch in~\cite{colton1996simple}. In the LSM, the values of an indicator function are used to identify points in the interior of the obstacle. Since the reconstruction of the domain is obtained by direct calculation of this indicator function, this method is classified as a sampling method or direct imaging method~\cite{colton1998inverse}. 

The linear sampling approach relies on finding level sets of an indicator function, and hence requires a crude estimate of the support of the obstacle for tabulating the indicator function.
In this work, we assume that the object is contained in the square $[-3,3]^2$, while all the domains are in fact contained in the square $[-1.5,1.5]^2$.

The original LSM was developed for the far-field pattern of scattered waves. However, in this work the scattered field is measured at distant receivers. So, we make the following adaptation for our setup. 
In a slight abuse of notation, let $\uscat(x,\theta)$ denote the scattered field at $x$ generated by the scattering of the incident plane wave $\uinc(x) = \exp{(i k x \cdot (\cos{(\theta), \sin{(\theta)}}))}$, and let $F$ denote the operator given by
\begin{equation}
F[g](x) = \int_{0}^{2\pi} \uscat(x,\theta) g(\theta) d\theta \, ,
\end{equation}
\revisea{where $g$ is known as the Herglotz wave function~\cite[Definition 3.26]{colton1998inverse}}. In particular, $F[g](x)$ is the solution to the Helmholtz equation with Dirichlet boundary conditions, and an incident field given by
$$
\uinc(x) = \int_{0}^{2\pi} \exp{(ikx \cdot (\cos{(\theta)}, \sin{(\theta)}))} g(\theta) \, d\theta \, .
$$
Assuming that measurements are made for sufficiently many equispaced angles in frequency, the operator $F[g](x)$ can be approximated to high accuracy via the trapezoidal rule as
\begin{equation}
F[g](x) \approx \frac{2\pi}{N_{d}} \sum_{\ell=1}^{N_{d}} \uscat(x,\theta_{\ell}) g(\theta_{\ell}) \, .
\end{equation} 
Let $\Phi(x,y)$ denote the scaled Helmholtz Green's function with wavenumber $k$ given by
\begin{equation*}
\Phi(x,y) = \exp(i\pi/4)\sqrt{\frac{\pi k}{2}}H_0^1(k|x-y|) \,,
\end{equation*}
where $H_{0}^{1}$ is the Hankel function of the first kind of order zero.
For each $x \in \mathbb{R}^{2}$, let 
$\boldsymbol{g}_{x} = [g_{x}(\theta_{1})\,; g_{x}(\theta_{2})\,; \ldots g_{x}(\theta_{N_{d}})]$ 
denote the discretized Herglotz wave function $g_{x}(\theta)$  satisfying 
\begin{equation}
F[g_{x}](x_{j}) = \Phi(x,x_{j}) \,, \quad j = 1,2,\ldots N_{t} \, ,
\end{equation}
where $x_{j}$ as before are the location of the receivers. 
If we reshape the measured data $A = \frac{2\pi}{N_{d}} \umeas$ as a $\mathbb{C}^{N_{t} \times N_{d}}$ matrix, then $\boldsymbol{g}_{x}$ is the solution of the following least squares problem 
\begin{equation}\label{eq:matrix_eqn}
\frac{2\pi}{N_{d}}   
\begin{bmatrix}
\umeas_{d_{1}}(x_{1}) & \umeas_{d_{2}}(x_{1}) & \ldots & \umeas_{d_{N_{d}}} (x_{1}) \\
\umeas_{d_{1}}(x_{2}) & \umeas_{d_{2}}(x_{2}) & \ldots & \umeas_{d_{N_{d}}} (x_{2}) \\
\vdots & \vdots & \ddots & \vdots \\
\umeas_{d_{1}}(x_{N_{t}}) & \umeas_{d_{2}}(x_{2}) & \ldots & \umeas_{d_{N_{d}}} (x_{N_{t}})
\end{bmatrix}
\cdot 
\begin{bmatrix}
g_{x} (\theta_{1}) \\
g_{x} (\theta_{2}) \\
\vdots \\
g_{x} (\theta_{N_{d}}) \\
\end{bmatrix}
=
\begin{bmatrix}
\Phi(x,x_{1}) \\ 
\Phi(x,x_{2}) \\
\vdots \\
\Phi(x,x_{N_{t}}) 
\end{bmatrix} \, .
\end{equation}
The above system of equations could be over- or under-determined based on the relation between $N_{t}$, and $N_{d}$. However, as long as $N_{t}$, and $N_{d}$ are large enough to resolve the scattered field measurements, the matrix tends to be exponentially ill-conditioned, and in practice, the above problem is solved in the least squares sense with Tikhonov regularization, i.e., 
\begin{equation}
\label{eq:lsm-lsq}
\boldsymbol{g}_{x} = \min_{\boldsymbol{g}} \| A\cdot \boldsymbol{g} - \boldsymbol{\Phi}_{x} \|^2  + \alpha^2 \| \boldsymbol{g} \|^2 \, , 
\end{equation}
where $\alpha$ is a regularization parameter, and $\boldsymbol{\Phi}_{x} = [\Phi(x,x_{1})\, ; \Phi(x,x_{2})\,; \ldots \Phi(x,x_{N_{t}})]$.
An indicator function denoted by $h(x)$ is then given by
\begin{equation}
h(x) = \log{ \left(\sqrt{\sum_{\ell=1}^{N_{d}} |g_{x}(\theta_{\ell})|^2} \right)} \, .
\end{equation}

Given the indicator function $h(x)$, the boundary is then defined via the set $h(x) = C$ for a some constant $C$. As noted in~\cite{cakoni2005qualitative}, 
$C$ is typically chosen using heuristics, and in this work we use the following heuristic. Suppose that the indicator function $h(x)$ is tabulated on a $200 \times 200$ equispaced grid of targets $x$ on $[-3,3]^2$. \revisea{Let $\text{dist}(S_{1},S_{2})$ denote the Chamfer distance \cite{barrow1977parametric} between two clouds of points $S_{1}$ and $S_{2}$, defined by}
\begin{equation}\label{eq:Chamfer}
\text{dist}(S_1, S_2) = \dfrac{1}{2|S_1|} \sum_{x \in S_1} \min_{y \in S_2} \norm{x-y} + \dfrac{1}{2|S_2|} \sum_{y \in S_2} \min_{x \in S_1} \norm{x-y}.
\end{equation}
Let $S_{j}$ represent a discretization of the level-set for $C_{j} = 7-0.2j$,  $j=0,1,2\ldots 15$. Then, $C=C_{j_{\textrm{opt}}}$ where $j_{\textrm{opt}}$ is the smallest $j$ such that $|\text{dist}(S_{j},S_{j+1}) - \text{dist}(S_{j+1},S_{j+2})|>0.1$. We found this heuristic to work best in practice for the range of frequencies, the family of star-shaped domains, and the particular choice of the level-set function $h(x)$ considered in this work.

The level set $S_{j_{\textrm{opt}}}$ obtained using this approach need not necessarily correspond to a star-shaped domain. We then find the set of coefficients $\bc$ such that the boundary $r(t;\bc)(\cos{(t)},\sin{(t}))$ is closest to the discretized level set $S_{j_{\textrm{opt}}}$ in the least squares sense as follows. Let $t_{\ell} \in S_{j_{\textrm{opt}}}$, $\ell=1,2,\ldots N$ denote the points on the level set. Let $\theta_{\ell} = \textrm{Arg}(t_{\ell}) \in (-\pi,\pi)$. Then the final reconstruction using the LSM method denoted by $\bclsm$ is given by
\begin{equation}
\bclsm = \min_{\bc} \sum_{\ell=1}^{N} ||t_{\ell}| - r(\theta_{\ell}; \bc)|^2 \, ,
\end{equation}
where $r(t;\bc)$ is the radius function given by~\cref{eq:rdef}.
The reconstruction of $\bclsm$ is not sensitive to the choice of the regularization parameter $\alpha$, as long as $\alpha$ is sufficiently small, but greater than $10^{-6}$, so that the least squares problem in~\cref{eq:lsm-lsq} can be solved in a stable manner in finite precision arithmetic. In all of the examples, we use $\alpha = 10^{-4}$.

\section{Neural network warm-start method} \label{sec:numerical_method}
In this section, we provide a detailed description of the neural network warm-start approach. In particular, we use a trained neural network to approximate the inverse obstacle scattering operator and obtain a guess for the coefficients
of the star-shaped obstacle $\bcnn$ such that $\cF_{k}(\bcnn)$ best fits the given scattered field measurements, i.e., $\bcnn$ tries to approximate the solution to~\cref{eq:optim_star}.
 The solution $\bcnn$ is subsequently refined using the Gauss-Newton iteration in~\Cref{alg}. In the following, we discuss the data generation for training the neural network, and the specific architecture used to approximate the inverse obstacle scattering operator. The code for data generation, model training, and warm-start method is available on the GitHub repository (see \url{https://github.com/frankhan91/learn-invscattering2d-obstacle}).

\subsection{Data generation}
\label{sec:data}
The neural network is trained on scattered data $\cF_{k}(\bc)$ with $\bc$ sampled from a distribution. Specifically, we let $c_{0} \sim \text{Unif}([1,1.2])$\footnote{$X\sim \text{Unif}([a,b])$ denotes a random variable uniformly distributed in the interval $[a,b]$.}, and $(c_{j}, c_{j+M}) =  (r\cos(\theta), r\sin(\theta))$,  for $j=1,2,\ldots M$,  where $r\sim \text{Unif}([0,0.1])$, and
$\theta \sim \text{Unif}([0,2 \pi])$. We omit any set of coefficients that result in a self-intersecting curve. 
In a slight abuse of notation, we reshape the vector of scattered data $\cF_{k}(\bc) \in \mathbb{C}^{N_{t} \cdot N_{d}}$ as a complex $N_{t} \times N_{d}$ matrix where each row corresponds to the scattered field at a fixed receiver for all incident directions, and each column corresponds to the scattered field at all receivers due to a fixed incident direction.

Unless stated otherwise, we assume that the receivers are equally spaced on the boundary of a circle (of radius 10) that is sufficiently far away from the obstacle, and that the scattered field is measured for equally spaced incident directions (the measurements could be full aperture or partial aperture). Moreover, we also assume that both the number of incident directions, and the number of receivers are sufficiently large so as to resolve the scattered field measurements as a function of the receiver angle on the circle, and of the angle of incidence.
For an aperture that is constant in frequency, this requirement translates to having scattered field measurements for $O(k)$ incident directions, and $O(k)$ receivers. 
Obtaining such highly-resolved scattered data tends to be difficult in practice. However, this idealized environment makes it feasible to study the behavior of the neural network warm-start approach as a function of the frequency of the incoming data, and the complexity of the shape of the obstacle without a simultaneous concern about the resolution of the scattered field measurements. 

Existing off-the-shelf tools for building and training neural networks work predominantly with real data. Thus, the real and the imaginary part of the scattered field data need to be dealt with separately. We address this issue by using only the real part of the scattered data. We further rescale and recenter the data in order to improve the training performance. In particular, let $\mu$ and $\sigma_{0}$ denote the element-wise mean and standard deviation of $\text{Re}(\cF_{k}(\bc))$. 
Then the input to the neural network is given by $(\text{Re}(\cF_{k})(\bc) - \mu)/\sigma_{0}$. 

\begin{remark}
We implemented a neural network where two input channels were used in the first layer corresponding to the real part and the imaginary part of the scattered data and found no significant improvement in performance. \revisea{This finding suggests that, in the problems considered in this study, the real part of the scattered data retains sufficient phase information for solving the inverse problem approximately. Thus, we present results using only the real part of the scattered data for brevity. Future research considering complex-valued neural networks, such as those discussed in \cite{Tygert2016,trabelsi2018deep}, may offer enhanced performance by incorporating the full complex measurement data when available.}
\end{remark}

\subsection{Network architecture}
One of the advantages of recasting the scattered data as a matrix and using the real part, is that the input to the neural network has the same input data format as an image. Thus, we can use a convolutional neural network that has proven to be successful in many image processing tasks. Broadly speaking, the network comprises of $L$ convolutional layers, with each of them consisting of a convolution, a non-linear activation, and a pooling operation, followed by $L'$ fully connected hidden layers. Finally, the coefficients $\bcnn$ are given by an affine transformation of the output of the last hidden layer.

Before presenting the architecture, we first review a few standard operations used in the network. Let $X \in \RR^{m\times n}$ and $W \in \RR^{m' \times n'}$ with $m' \le m$ and $n' \le n$. Then the cross-correlation of $W$ and $X$, denoted by $W \star X\in \RR^{(m - m' +1) \times (n - n' +1)}$, is given by
$$(W \star X)_{i,j} = \sum_{i'=1}^{m'} \sum_{j'=1}^{n'} W_{i', j'}  X_{(i'+i-1), (j'+j-1)}, \quad 1 \leq i \leq m-m'+1,\;1\leq j \leq n -n' + 1.$$
Let $\mathrm{pad}(X,p) \in \RR^{(m+2p) \times (n+2p)}$ denote the padding operator with padding size $p$ given by
$$\parentheses{ \mathrm{pad}(X,p) }_{(p+1):(p+m),(p+1):(p+n)} = X.$$
For full aperture measurements, a periodic version of the padding operator is used since the input data is periodic both in the angle of incidence, and also in the location of the receivers, i.e.,
$$\parentheses{ \mathrm{pad}(X,p) }_{i,j} = \parentheses{ \mathrm{pad}(X,p) }_{i+m,j+n} \, , \quad 1\leq i,j \leq 2p \,. $$
On the other hand, for partial aperture measurements, we use zero padding instead of periodic padding. In a slight abuse of notation, let ``$\textrm{pad}$'' denote the zero padding operator as well, which is given by
$$\parentheses{ \mathrm{pad}(X,p) }_{i,j} = 0 \, , \quad 1\leq i,j \leq p \, , \quad m+p+1 \leq i \leq m+2p \, , \quad n+p+1 \leq j \leq n+2p \,. $$
Let $\mathrm{pool}(X) \in \RR^{\lfloor m/2 \rfloor \times \lfloor n/2\rfloor}$ denote the average pooling operator given by
\begin{equation}
\mathrm{pool}(X)_{i,j} = \frac{1}{4} \left( X_{2i-1,2j-1} + X_{2i-1,2j} + X_{2i,2j-1} + X_{2i,2j} \right).
\end{equation}
Finally, let ReLU$(x)$ denote the standard rectified linear unit given by $\mathrm{ReLU}(x)=\max(x,0)$ in the componentwise sense.

Given the above definitions, the specifics of each of the layers are discussed below. The neural network is parametrized via the following $6$ sets of parameters: the weight matrices and biases in the  convolutional layers, the weight matrices and biases in the fully connected layers, and finally the weight matrix and bias to obtain the output from the last hidden layer.

Suppose that each convolutional layer has $N_{c}$ channels in the output, and let $X^{\ell}_{k}$ denote the output of the $\ell$th convolution layer and the $k$th channel. 
With this notation, the input data is $X^{0}_{1} = (\text{Re}(\cF_{k})(\bc) - \mu)/\sigma_{0}$, \reviseb{recalling that $\mu$ and $\sigma_{0}$ denote the element-wise mean and standard deviation of $\text{Re}(\cF_{k}(\bc))$}.
Let $W^{\ell}_{j,k} \in \RR^{n_{K} \times n_{K}}$, and $b_{k}^{\ell} \in \RR$, $\ell=1,\dots, L$, $j=1,\ldots N_{c}$, and $k=1,\ldots N_{c}$ denote the \reviseb{trainable} weight matrices  and biases in the convolutional layer respectively ($\ell$ indexes the layer number, $j$ and $k$ index the channel, and $n_{K}$ is the kernel width).
The output of the $\ell$th convolution layer with inputs $\{ X_{j}^{\ell-1} \}_{j=1}^{N_{c}}$ (with the understanding that $N_{c} = 1$, if $\ell=1$) is given by
\begin{equation}
X_{k}^{\ell} = \textrm{pool} \left( \textrm{ReLU} \left( \sum_{j=1}^{N_{c}} W_{j,k}^{\ell} \star \textrm{pad} \left( X_{j}^{\ell-1},p \right) \oplus b_{k}^{\ell}  \right) \right) \, .
\end{equation}
Here $\oplus$ denotes componentwise addition: if $X_{j}^{\ell-1} \in \RR^{m\times n}$ for all $j$, then $\sum_{j=1}^{N_{c}} W_{j,k}^{\ell} \star \textrm{pad} \left( X_{j}^{\ell-1},p \right) \in \RR^{(m+2p-n_{K}+1) \times (n+2p-n_{K}+1)}$, and its each element is added by the scalar $b_{k}^{\ell}$.
In the remainder, we assume that $N_{t}$, $N_{d}$ divide $2^{L}$, and set the kernel width $n_{K} = 2p + 1$, which results in $X_{k}^{\ell} \in \RR^{N_{t}/2^{\ell} \times N_{d}/2^{\ell}}$ for all $k,\ell$.

\begin{figure*}[!htb]
    \centering
    \includegraphics[width=0.99\textwidth]{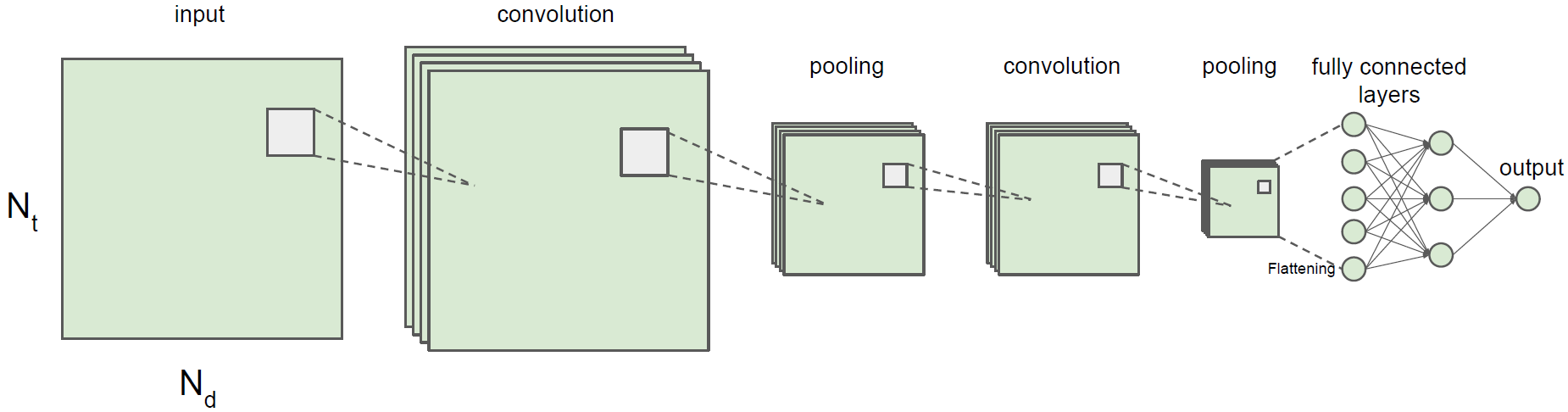}
    \caption{The structure of the neural network: $2$ layers of convolution-ReLu-pooling, followed by $2$ fully connected hidden layers.}
    \label{fig:nn}
\end{figure*}

After obtaining the result from the convolutional layers, we flatten it to a single vector denoted by $x^{0} \in \RR^{N_{t} \cdot N_{d} \cdot N_{c}/4^{L}}$. Let $W^{\ell} \in \RR^{N_{\ell} \times N_{\ell-1}}$, and $b^{\ell} \in \RR^{N_{\ell}}$, $\ell \in 1,2 \ldots L'$, denote the \reviseb{trainable} weight matrices and biases in the fully-connected layers with $N_{0} = N_{t} \cdot N_{d} \cdot N_{c}/4^{L}$.
The output of the $\ell$th fully-connected layer is then given by
\begin{equation}
x^{\ell} = \textrm{ReLU} \left( W^{\ell} x^{\ell-1} + b^{\ell} \right) \, .
\end{equation}
The final coefficients $\bcnn$ are then given by
\begin{equation}
\bcnn = \widetilde{W} x^{L'} + \tilde{b} \, ,
\end{equation}
with the \reviseb{trainable weight matrix $\widetilde{W} \in \RR^{N_{L'} \times (2M+1)}$, and the bias $\tilde{b} \in \RR^{2M+1}$ in the last linear layer}.
This architecture with $2$ convolutional layers, and $2$ fully-connected layers with $M=0$ is illustrated in Figure \ref{fig:nn}.

The neural network is trained on the normalized dataset generated using the stochastic gradient descent (SGD) method with momentum. 
Let $\Nepo$ denote the number of epochs for which the network is trained. 
In all the experiments, we set the batch size to $100$; the learning rate is set to $0.16$ for the first $(\Nepo - 100)$ epochs, and to $0.08$ for the last $100$ epochs. \revisea{The neural network and the associated training process are implemented utilizing the open-source machine learning framework PyTorch~\cite{NEURIPS2019_9015}.}

\section{Numerical results}\label{sec:numerical_results}
The performance of solvers for the solution to~\cref{eq:optim_star} is sensitive to the frequency of the incident waves $k$, and the complexity of the radius function defining the star-shaped obstacles. In this section, we explore the impact of $k$ on the following $5$ different approaches.
\begin{itemize}
\item \textbf{GN:} The solution is computed using the Gauss-Newton approach in~\Cref{alg} initialized with a unit disk centered at the origin. 
\item\textbf{LSM prediction:} The solution is the output of the linear sampling method $\bclsm$ described in~\Cref{sec:LSM}. 
 \item\textbf{LSM refined:} The solution is computed using the Gauss-Newton approach in~\Cref{alg} initialized with $\bc^{(0)} = \bclsm$.
\item\textbf{DL prediction:} The solution is given by the output of the neural network $\bcnn$ described in~\Cref{sec:numerical_method} (without any refinement using the Gauss-Newton approach).
  \item\textbf{DL refined:} The solution is computed using the Gauss-Newton approach in~\Cref{alg} initialized with $\bc^{(0)} = \bcnn$.
 \end{itemize}

\begin{table}[h!]\centering 
\begin{tabular}{c|ccccccc}
\toprule
$M$ & ($N_{d}$, $N_{t}$) & $\Ntrain$ & $\Nepo$ & $N_{c}$ & $(p,  n_{K})$ & $N_{1}$ & $N_{2}$\\
\midrule
$5$ & (48,48) & $500$  & $1000$ & $5$ & $(2,5)$ & $50M$ & $10M$\\
\midrule 
$10$ & (48,48) & $2000$  & $1000$ & $10$ & $(2,5)$ & $50M$ & $10M$\\
\midrule
$20$ & (100,100) & $80000$  & $2000$ & $20$ & $(4,9)$ & $50M$ & $10M$\\
\bottomrule
\end{tabular}
\caption{Parameters defining the network architecture, data set size, and network training. $M$ is the maximum Fourier content of the radius function defining the boundary, $\Nepo$ denotes the number of SGD iterations for training the network, $\Ntrain$ denotes the number of training data used to train the neural network. $N_{c}$ is the number of channels in the convolutional layer, $p$ is the padding width, $n_{K}$ is the kernel width of the weight matrices in the convolutional layers, and $N_{1}$, and $N_{2}$ define the weight matrices in the fully-connected layers. \reviseb{The experiments were run on a NVIDIA A100 Tensor Core GPU. The training time for the three models are $22.8s$, $67.1s$, and $21004.2s$.}}
\label{tab:para}
\end{table}

The detailed parameters for the numerical implementation of different methods are summarized below.
We run the Gauss-Newton approach for a maximum of $N_{\textrm{max}} = 20$ iterations, 
and stop if the update tolerance is below $\ve_{s} = 5 \times 10^{-8}$, or if the residue tolerance is below
$\ve_{r} = \num{1e-6}$.
For both the GN approach discussed above, and the computation of 
$\bclsm$, we use $N_{t} = 200$ receiver locations on the circle of radius $10$, 
and $N_{d} = 200$ incident waves. 
This is more than double the maximum number of receiver locations, and incident waves used for the DL based methods, since neither LSM or GN use any training data for computing the solution.
For the neural network, we have two convolutional layers ($L=2$) and two fully connected hidden layers ($L'=2$). The other parameters of the neural network architecture and training are summarized in Table~\ref{tab:para}, which depend on the maximum Fourier content $M$ of the radius function defining the star-shaped domains.  

In the following, the error in reconstruction is measured via the relative error in terms of the $\ell^2$ distance in $\RR^{2M+1}$ between the computed Fourier coefficients vector and the Fourier coefficients of the obstacle corresponding to the measured data. The qualitative behavior of error in reconstruction was identical for other monitor functions such as the Chamfer distance discussed in~\Cref{sec:LSM}. The numerical results comparing all five approaches are summarized in~\Cref{tab:summary}. The results are computed by averaging the relative errors on 50 test geometries sampled from the same distribution of geometries for generating the training data detailed in Section~\ref{sec:data}. 
\reviseb{The average running time per obstacle geometry for all five approaches is reported in~\Cref{tab:time}.}

\begin{table}[h!]\centering 
{\footnotesize
\begin{tabular}{l@{\hskip0pt}l@{\hskip0pt}l|ccccc}
\toprule
& & & GN & LSM prediction & LSM refined & DL prediction & DL refined\\
\midrule
$M=5$, & $k=5$ & & $10.78\%$ & $2.17\%$ & $0.40\%$ & $5.40\%$ & $0$\\

$M=10$,~ & $k=10$\phantom{,}~ & & $21.29\%$ & $13.04\%$ & $7.72\%$ & $3.57\%$ & $0$\\

$M=20$, & $k=30$ & & $24.53\%$ & $32.21\%$ & $32.31\%$ & $4.30\%$ & $1.32\%$\\
\midrule 

$M=5$, & $k=5$, & $\sigma=0.05$ & $10.78\%$ & $8.72\%$ & $4.33\%$ & $5.40\%$ & $0.083\%$\\

$M=10$, & $k=10$, & $\sigma=0.05$ & $21.29\%$ & $22.10\%$ & $21.30\%$ & $3.57\%$ & $0.054\%$\\

$M=20$, & $k=30$, &$\sigma=0.05$ & $24.54\%$ & $34.41\%$ & $34.57\%$ & $4.30\%$ & $1.35\%$\\

$M=20$, & $k=30$, &$\sigma=0.15$ & $24.55\%$ & $49.54\%$ & $49.42\%$ & $4.30\%$ & $1.39\%$\\

\midrule 
$M=10$, & $k=10$, & partial data & $19.46\%$ & $21.37\%$ & $18.38\%$ & $6.02\%$ & $0$\\
\bottomrule
\end{tabular}
}
\caption{Summary of the numerical results in terms of average relative errors over 50 test geometries. For the GN, LSM prediction, and LSM refined approach, the scattered data was measured at $N_{t} = 200$ receivers, and for $N_{d} = 200$ incident directions, while for the DL based approaches, the number of receivers and incident directions were chosen based on the parameters in Table~\ref{tab:para}. The $0$s in the table indicate that the algorithm reaches machine precision.}
\label{tab:summary}
\end{table}

\begin{table}[h!]\centering 
{\footnotesize
\begin{tabular}{l@{\hskip0pt}l@{\hskip0pt}l|ccccc}
\toprule
& & & GN & LSM prediction & LSM refined & DL prediction & DL refined\\
\midrule
$M=5$, & $k=5$ & & $9.1$ & $21.0$ & $23.7$ & $2.2$ & $4.3$\\

$M=10$,~ & $k=10$\phantom{,}~ & & $29.9$ & $20.3$ & $42.8$ & $3.4$ & $10.3$\\

$M=20$, & $k=30$ & & $83.7$ & $96.4$ & $186.6$ & $4.5$ & $50.2$\\
\midrule 

$M=5$, & $k=5$, & $\sigma=0.05$ & $9.2$ & $19.6$ & $24.3$ & $2.3$ & $4.7$\\

$M=10$, & $k=10$, & $\sigma=0.05$ & $28.9$ & $18.2$ & $46.4$ & $3.7$ & $10.7$\\

$M=20$, & $k=30$, &$\sigma=0.05$ & $78.0$ & $85.0$ & $164.3$ & $4.9$ & $48.1$\\

$M=20$, & $k=30$, &$\sigma=0.15$ & $78.2$ & $138.2$ & $257.9$ & $4.3$ & $48.6$\\

\midrule 
$M=10$, & $k=10$, & partial data & $14.6$ & $16.7$ & $28.7$ & $0.6$ & $7.1$\\
\bottomrule
\end{tabular}
}
\caption{\reviseb{Average running time (in seconds) per obstacle geometry for all five approaches. All the experiments were conducted on an Intel Xeon Gold 6148 CPU (2.40GHz). The DL prediction time does not include the training time, while the refined methods include prediction time.}}
\label{tab:time}
\end{table}

We remark that the relative errors of some approaches in some cases have a bimodal distribution. For example, when $M=5$, $k=5$ (the first row in~\Cref{tab:summary}), approximately half solutions on the 50 testing cases, provided by the Gauss-Newton approach, have zero relative error while another half of the solutions have relative errors around $20\%$. To better present the statistics of errors, we also report in~\Cref{tab:percentage} the percentages of cases whose relative error is lower than $1\%$ for three approaches (GN, LSM refined, DL refined). 

\begin{table}[h!]\centering 
{\footnotesize
\begin{tabular}{l@{\hskip0pt}l@{\hskip0pt}l|ccc}
\toprule
& & & GN  & LSM refined & DL refined\\
\midrule
$M=5$, & $k=5$ & & $44\%$ & $98\%$ & $100\%$\\

$M=10$,~ & $k=10$\phantom{,}~ & & $12\%$ & $52\%$ & $100\%$\\

$M=20$, & $k=30$ & & $0$ & $0$ & $66\%$\\
\midrule 

$M=5$, & $k=5$, & $\sigma=0.05$ & $44\%$ & $86\%$ & $100\%$\\

$M=10$, & $k=10$, & $\sigma=0.05$ & $12\%$ & $22\%$ & $100\%$\\

$M=20$, & $k=30$, &$\sigma=0.05$ & $0$ & $0$ & $64\%$\\

$M=20$, & $k=30$, &$\sigma=0.15$ & $0$ & $0$ & $64\%$\\

\midrule 
$M=10$, & $k=10$, & partial data & $40\%$ & $38\%$ & $100\%$\\
\bottomrule
\end{tabular}
}
\caption{Percentage of cases whose relative error is lower than $1\%$, following the setting described in~\Cref{tab:summary}.}
\label{tab:percentage}
\end{table}

\subsection{Full-aperture noiseless data}
We first compare the $5$ approaches for full-aperture and noiseless data, for $(k,M) = (5,5), (10,10)$, and $(30,20)$. For each of these configurations, $\umeas_{k} = \cF(\bcmeas)$, where $\bcmeas$ is generated from the same distribution as the training distribution: $c_{0} \sim \text{Unif}([1,1.2])$, and $(c_{j}, c_{j+M}) =  (r\cos(\theta), r\sin(\theta))$,  for $j=1,2,\ldots M$,  where $r\sim \text{Unif}([0,0.1])$, $\theta\sim \text{Unif}([0,2\pi])$. 

\begin{figure*}[!htb]
    \centering
    \includegraphics[width=\textwidth]{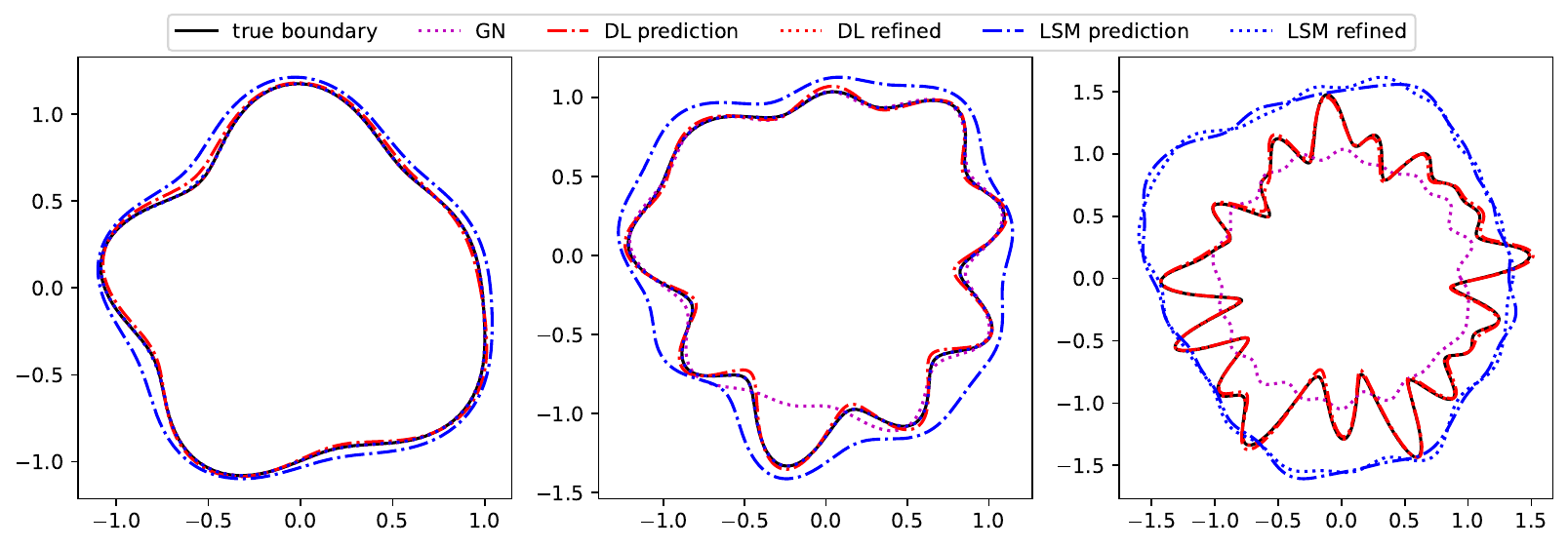}
    \caption{{\bf Noiseless full-aperture measurements}: Reconstructions obtained using the $5$ approaches: GN, LSM prediction, LSM refined, DL prediction, and DL refined for full-aperture noiseless data. (Left) $(k,M) = (5,5)$, (middle) $(k,M) = (10,10)$, (right) $(k,M) = (30,20)$. 
Scattered field measurements were made at $N_{t} = 200$ receivers, for incident waves from $N_{d} = 200$ directions for the GN, LSM prediction, and LSM refined approaches, while for the DL based approaches $N_{t}$, and $N_{d}$ were chosen on a frequency dependent manner: $N_{t} = N_{d} = 48$, for $(k,M) = (5,5)$, and $(10,10)$, and $N_{t} = N_{d} = 100$, for $(k,M) = (30,20)$.}
    \label{fig:first}
\end{figure*}

Referring to~\Cref{tab:summary}, for the lowest frequency configuration $(k,M) = (5,5)$, $\bclsm$, and $\bcnn$ are sufficiently close to the true solution with an error of $2.17 \%$, and $5.40\%$ respectively.
Gauss-Newton running with the unit circle converge to the exact solution in about half testing cases, while Gauss-Newton running with both $\bclsm$ and $\bcnn$ as the initial guess converge to the true solution accurately, with only one exception for the LSM method.
As the frequency is increased to $k=10$, the GN approach does not converge to the correct solution, while $\bclsm$ has a $13.04\%$ error but still lies in the local set of convexity of the true solution $\bcmeas$ in about half cases ($52\%$ reported in~\Cref{tab:percentage}). Thus, when we refine the solution using the Gauss-Newton iteration with $\bclsm$ as the initial guess, it is still possible to get a perfect recovery (shown in the middle of ~\Cref{fig:first}). The DL prediction $\bcnn$ provides a much better initial guess than $\bclsm$ as is evidenced by the corresponding errors, and the solution upon refinement converges to $\bcmeas$ as well.
Finally, all methods, except the DL prediction $\bcnn$ and the DL refined solution, worsen significantly when the frequency is increased further. The solution obtained via GN, LSM prediction or LSM refined are $O(1)$ incorrect.
The lack of convergence with increasing frequency is further illustrated through Figure~\ref{fig:first}, where we plot the reconstructions for $(k,M)= (5,5), (10,10)$, and $(30,20)$.

\subsection{Effect of measurement noise}
In this section, we study the impact of measurement noise on the performance of all $5$ approaches. 
In particular, let 
$$
\umeas_{k,d_{\ell}}(x_{j}) = \left(\cF_{k,d_{\ell}}(\bcmeas) \right)_{j} \left( 1+ \xi_{l,j} \cdot \exp\parentheses{{i\chi_{l,j}}} \right),\quad  \xi_{l,j} \sim \mathrm{Unif} ([0,2\sigma]), \chi_{l,j}\sim \mathrm{Unif} ([0,2\pi]),$$
$\ell = 1,2,\ldots N_{d}$, $j=1,2,\ldots N_{t}$,
i.e., we consider multiplicative noise with average amplitude $\sigma$, independently applied to each measurement.
Note that even though the measurement data was noisy, the neural network was trained using noiseless measurements. 

\begin{figure*}[h!]
    \centering
    \includegraphics[width=\textwidth]{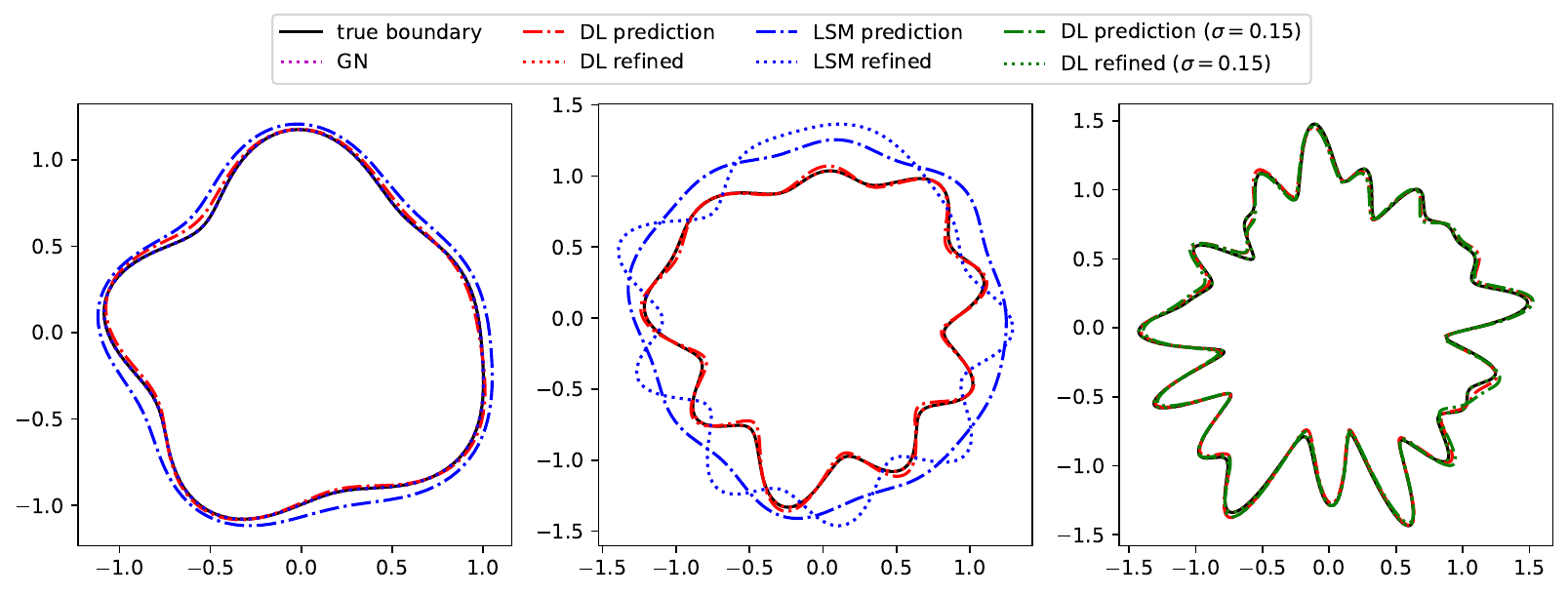}
    \caption{ {\bf Results with noisy measurements}: Reconstructions with noisy measurements using all $5$ approaches with mean multiplicative noise of $\sigma = 0.5$  for (left) $(k,M) = (5,5)$, and middle $(k,M) = (10,10)$. On the right, we plot the results for the DL based approaches for $(k,M) = (30,20)$, and $\sigma = 5\%$ and $\sigma = 15\%$.}
    \label{fig:noise}
\end{figure*}

The performance of the various algorithms with a multiplicative noise of $5\%$ is similar to the case of noiseless measurements. For low frequency $(k,M) = (5,5)$, all of the algorithms recover the shape of the obstacle to high accuracy, similar to the situations without noise, and as the frequency is increased, the DL based approaches result in a robust reconstruction of the obstacle, while the GN, and LSM based approaches fail to converge. The DL based approaches recover the shape of the obstacle to high fidelity even when $\sigma = 15\%$. The reconstructions obtained using all $5$ approaches are plotted in Figure~\ref{fig:noise}.

\begin{remark}
Referring to Table~\ref{tab:summary}, it is interesting to note that the error in the coefficients obtained using the Gauss-Newton approach with and without noise are the same to $4$ significant digits, This is evident from the plots of the two reconstructions as well. However, the number of iterations in the Gauss-Newton iteration and final iterations were different between the runs. The similarity in the final reconstruction could in part be explained by the smallness of the noise, and averaging happening across various measurements. 
\end{remark}

\subsection{Partial aperture data}
\begin{figure*}[!htb]
    \centering
    \includegraphics[width=0.7\textwidth]{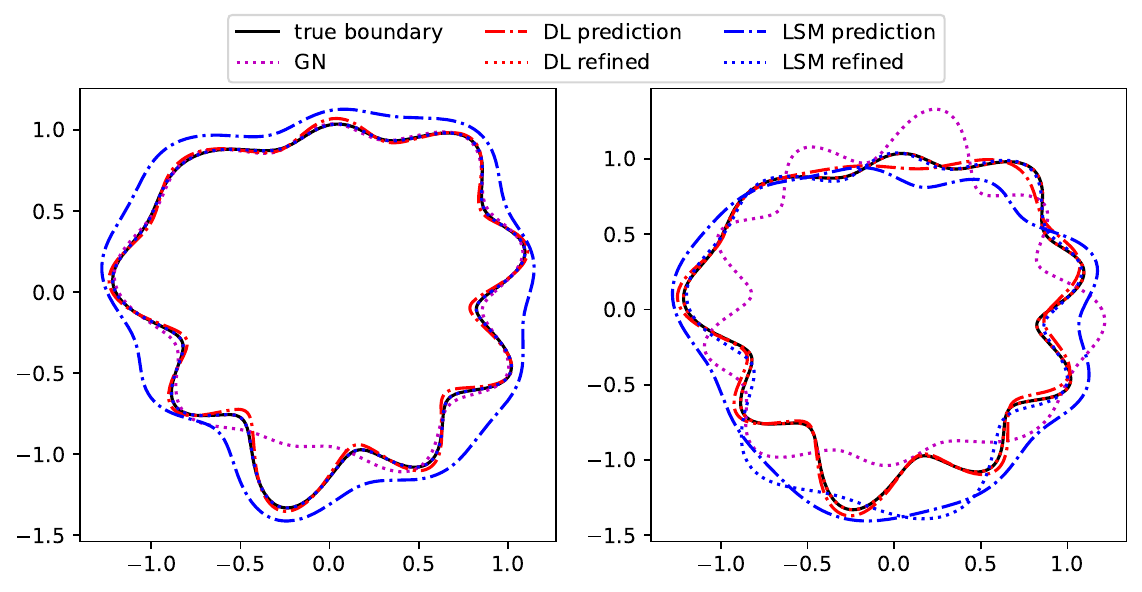}
    \caption{ {\bf Results with partial aperture data}: Reconstructions for $(k,M) = (10,10)$ for (left)  full-aperture noiseless measurements (same as Figure \ref{fig:first} middle), and (right) partial aperture noiseless measurements. For the plot on the right, measurements were made at receivers located at $r=10$, and equally spaced angles between $[0,\pi)$, and incident directions which were equally spaced in $[0,\pi)$.}
    \label{fig:partial}
\end{figure*}

In this section, we compare the performance of the algorithms with partial 
aperture data for $(k,M) = (10,10)$. The receivers 
are located on $x_{j} = 10(\cos{(j \pi /N_{t})}, \sin{(j \pi/N_{t})})$, with $j=1,2,\ldots, N_{t}$, and the incidence 
directions are $d_{\ell} = (\cos{(\ell \pi/N_{d})} , \sin{(\ell \pi/N_{d})})$, with
$\ell=1,2,\ldots N_{d}$. The neural network was 
trained with partial aperture data as well. 
While LSM with refinement converges when full aperture measurements are available, 
it fails to converge to the correct solution with partial aperture measurements. 
On the other hand, while the accuracy of the DL prediction method deteriorates slightly, 
it still lies in the basin of attraction of the true solution, and 
we recover the obstacle accurately using DL refined. In Figure~\ref{fig:partial},  
we plot the results for all five methods.

\subsection{Training data scaling with shape complexity}
\begin{figure*}[h!]
    \centering
    \includegraphics[width=0.95\textwidth]{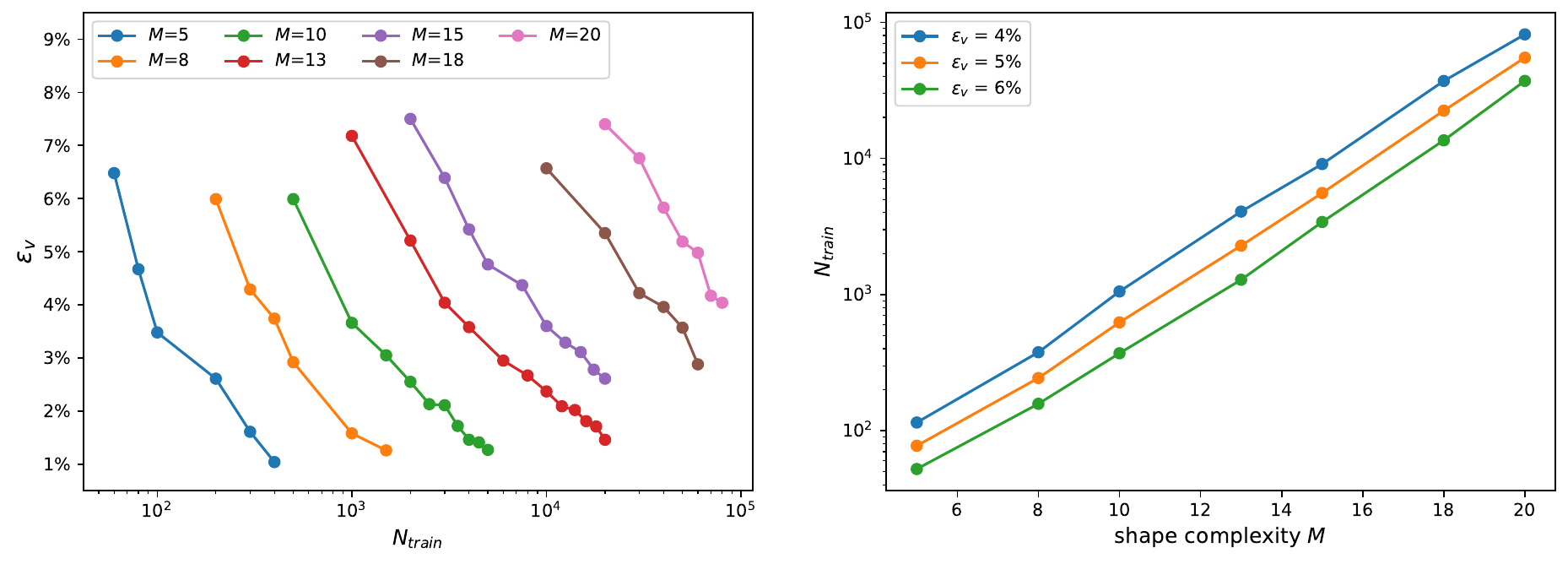}
    \caption{The error curves. Left: final validation errors for different $M$ and different number of data. Right: number of data required to achieve certain accuracy for different $M$. The slope of all the three lines is $0.192$ with a score larger than $0.997$.}
    \label{fig:error_curve}
\end{figure*}

Finally, we conclude the section with a study of the scaling of the number of 
training data $\Ntrain$ with shape complexity $M$ (and hence frequency $k$, since $M= O(k)$) in order to achieve a given validation error $\ve_{v}$.
The validation error $\ve_{v}$ is the averaged relative error 
$\|\bcnn - \bcmeas\|/\|\bcmeas\|$ for 500 sets of test coefficients $\bcmeas$. 
Let the maximum Fourier content of the radius function $r(t;\bc)$ be $M = k - 10$. We measure the scattered data at $N_{t} = 100$ receivers, and for $N_{d} = 100$ incident directions. The choice of $M=k-10$ has been made based on the success of the DL based approaches in this setting. The other network parameters are set as follows, $p=4$, $n_{K} = 9$, and the widths defining the weight matrices in the fully connected layers to $N_{1} = 50M$, and $N_{2} = 10M$. 
In Figure~\ref{fig:error_curve}, we plot the validation error with the increasing number 
of training samples as a function of $M$. The plot shows that in order to achieve a 
fixed validation error, 
$\Ntrain$ grows exponentially in frequency. 
This is further illustrated by computing a linear fit between $\log{\Ntrain}$ and $\ve_{\textrm{v}}$ for each $M$ and $\Ntrain$ as a function of $M$ for the fixed validation error. In particular, we note that $\Ntrain \approx 7.68 \times 1.556^{M}$ for $\ve_{v} = 5\%$.

\section{Conclusion}\label{sec:conclusion}
In this work, we presented a neural network warm-start approach for the solution of the inverse acoustic obstacle scattering problem restricted to star-shaped domains in two dimensions. 
A trained neural network is used to obtain a high fidelity guess for the shape of the obstacle from scattered field measurements at distant receivers, which is then further refined using a Gauss-Newton based iterative approach.

For noiseless full-aperture data, the solution obtained using the trained neural network was able to achieve an accuracy comparable to the validation error, which would then converge to the exact solution when refined using the Gauss-Newton approach. The algorithm remained robust with increasing frequency, where a tighter initial guess is required owing to the shrinking local set of convexity in the vicinity of the true solution. The qualitative behavior of solutions obtained using this approach was not impacted by the presence of noise in scattered field measurements, or if only partial aperture measurements were available. 
We compared the efficiency of this approach to three other approaches, the Gauss-Newton approach with the unit circle as an initial guess, a linear sampling method adapted to the problem setup, and a linear sampling warm-start approach, and found that in all the regimes considered in this work, the neural network warm-start approach outperformed all the other methods for sufficiently high frequency inverse problems.

However, the number of training samples required to train the network to achieve a fixed validation error scales exponentially with the frequency of the scattered field measurements and the complexity of the obstacles considered. 
Based on the complicated nature of the optimization landscape of the problem, this behavior is expected since the number of local minima around the true obstacle, also scales exponentially with frequency.
A natural question arises: are there other network architectures that can break this exponential scaling of complexity? If not, can the network architecture be improved to minimize the number of training samples required? In many experimental settings, the aperture of the scattered field measurements is severely restricted and often times only backscatter measurements from the obstacle are available. Traditional approaches like the Gauss-Newton approach and the linear sampling methods are typically not robust in this regime. The impact of having such limited aperture data on the robustness of the neural network warm-start approach both in terms of accuracy of recovery, and the number of training samples required to train the network is an area of ongoing research.

The approach presented in this paper naturally extends to obstacles with sound-hard, impedance or transmission boundary conditions as well. This approach could also be adapted for the recovery of volumetric properties such as variations in sound speed, or density of the object. This line of inquiry is also being vigorously pursued and will be reported at a later date. The ability to develop efficient neural network based approaches in two dimensions and understanding their training complexity with frequency will play a critical role in the extension of these methods to inverse acoustic and electromagnetic problems in three dimensions.

\section{Acknowledgments}
The work of C. Borges was supported in part by the Office of Naval Research under award number N00014-21-1-2389. We would like to thank Travis Askham, Leslie Greengard, and Jeremy Hoskins for many useful discussions.

\bibliographystyle{elsarticle-num} 
\bibliography{ref.bib}

\begin{thebibliography}{10}
\expandafter\ifx\csname url\endcsname\relax
  \def\url#1{\texttt{#1}}\fi
\expandafter\ifx\csname urlprefix\endcsname\relax\def\urlprefix{URL }\fi
\expandafter\ifx\csname href\endcsname\relax
  \def\href#1#2{#2} \def\path#1{#1}\fi

\bibitem{kuchment2014radon}
P.~Kuchment, \href{https://doi.org/10.1137/1.9781611973297}{{T}he {R}adon
  {T}ransform and {M}edical {I}maging}, CBMS-NSF Regional Conference Series in
  Applied Mathematics, Society for Industrial and Applied Mathematics, 2014.
\newline\urlprefix\url{https://doi.org/10.1137/1.9781611973297}

\bibitem{collins1995nondestructive}
R.~Collins,
  \href{https://books.google.com/books?id=E55RFcv-rGMC}{{N}ondestructive
  {T}esting of {M}aterials}, Studies in applied electromagnetics and mechanics,
  IOS Press, 1995.
\newline\urlprefix\url{https://books.google.com/books?id=E55RFcv-rGMC}

\bibitem{engl2012inverse}
H.~Engl, A.~K. Louis, W.~Rundell,
  \href{https://doi.org/10.1007/978-3-7091-6521-8}{{I}nverse {P}roblems in
  {M}edical {I}maging and {N}ondestructive {T}esting: {P}roceedings of the
  {C}onference in {O}berwolfach, {F}ederal {R}epublic of {G}ermany, {F}ebruary
  4--10, 1996}, Springer Vienna, 2012.
\newline\urlprefix\url{https://doi.org/10.1007/978-3-7091-6521-8}

\bibitem{Ustinov2014}
E.~Ustinov,
  \href{http://dx.doi.org/10.1007/978-0-387-36699-9_54}{{E}ncyclopedia of
  {R}emote {S}ensing}, Springer New York, New York, NY, 2014, Ch. Geophysical
  Retrieval, Inverse Problems in Remote Sensing, pp. 247--251.
\newblock \href {https://doi.org/10.1007/978-0-387-36699-9_54}
  {\path{doi:10.1007/978-0-387-36699-9_54}}.
\newline\urlprefix\url{http://dx.doi.org/10.1007/978-0-387-36699-9_54}

\bibitem{cheney2009fundamentals}
M.~Cheney, B.~Borden,
  \href{https://doi.org/10.1137/1.9780898719291}{{F}undamentals of {R}adar
  {I}maging}, CBMS-NSF Regional Conference Series in Applied Mathematics,
  Society for Industrial and Applied Mathematics, 2009.
\newline\urlprefix\url{https://doi.org/10.1137/1.9780898719291}

\bibitem{beilina2015globally}
L.~Beilina, N.~T. Thanh, M.~V. Klibanov, J.~B. Malmberg,
  \href{https://doi.org/10.1016/j.cam.2014.11.055}{Globally convergent and
  adaptive finite element methods in imaging of buried objects from
  experimental backscattering radar measurements}, Journal of Computational and
  Applied Mathematics 289 (2015) 371--391.
\newline\urlprefix\url{https://doi.org/10.1016/j.cam.2014.11.055}

\bibitem{thanh2015imaging}
N.~T. Thanh, L.~Beilina, M.~V. Klibanov, M.~A. Fiddy,
  \href{https://doi.org/10.1137/140972469}{Imaging of buried objects from
  experimental backscattering time-dependent measurements using a globally
  convergent inverse algorithm}, SIAM Journal on Imaging Sciences 8~(1) (2015)
  757--786.
\newline\urlprefix\url{https://doi.org/10.1137/140972469}

\bibitem{chen1997inverse}
Y.~Chen, \href{https://doi.org/10.1088/0266-5611/13/2/005}{Inverse scattering
  via heisenberg's uncertainty principle}, Inverse problems 13~(2) (1997) 253.
\newline\urlprefix\url{https://doi.org/10.1088/0266-5611/13/2/005}

\bibitem{borges2020inverse}
C.~Borges, J.~Lai, \href{https://doi.org/10.1088/1361-6420/abac9b}{Inverse
  scattering reconstruction of a three dimensional sound-soft axis-symmetric
  impenetrable object}, Inverse Problems 36~(10) (2020) 105005.
\newblock \href {https://doi.org/10.1088/1361-6420/abac9b}
  {\path{doi:10.1088/1361-6420/abac9b}}.
\newline\urlprefix\url{https://doi.org/10.1088/1361-6420/abac9b}

\bibitem{chen1995recursive}
Y.~Chen,
  \href{http://cpsc.yale.edu/sites/default/files/files/tr1088.pdf}{Recursive
  linearization for inverse scattering}, Mathematical and numerical aspects of
  wave propagation (Golden, CO, 1998) (1995) 114--117.
\newline\urlprefix\url{http://cpsc.yale.edu/sites/default/files/files/tr1088.pdf}

\bibitem{bao2005inverse}
G.~Bao, P.~Li, \href{https://doi.org/10.1137/040607435}{Inverse medium
  scattering problems for electromagnetic waves}, SIAM Journal on Applied
  Mathematics 65~(6) (2005) 2049--2066.
\newline\urlprefix\url{https://doi.org/10.1137/040607435}

\bibitem{bao2012shape}
G.~Bao, P.~Li, \href{https://doi.org/10.1515/9783110259056.283}{Shape
  reconstruction of inverse medium scattering for the Helmholtz equation},
  Vol.~56, Higher Education Press Beijing, 2012.
\newline\urlprefix\url{https://doi.org/10.1515/9783110259056.283}

\bibitem{sini2012inverse}
M.~Sini, N.~T. Thanh, \href{https://doi.org/10.3934/ipi.2012.6.749}{Inverse
  acoustic obstacle scattering problems using multifrequency measurements},
  Inverse Problems \& Imaging 6~(4) (2012) 749.
\newline\urlprefix\url{https://doi.org/10.3934/ipi.2012.6.749}

\bibitem{chaillat2012faims}
S.~Chaillat, G.~Biros,
  \href{https://doi.org/10.1016/j.jcp.2012.02.006}{Fa{IMS}: A fast algorithm
  for the inverse medium problem with multiple frequencies and multiple sources
  for the scalar helmholtz equation}, Journal of Computational Physics 231~(12)
  (2012) 4403--4421.
\newline\urlprefix\url{https://doi.org/10.1016/j.jcp.2012.02.006}

\bibitem{borges2015inverse}
C.~Borges, L.~Greengard, \href{https://doi.org/10.1137/140982787}{Inverse
  obstacle scattering in two dimensions with multiple frequency data and
  multiple angles of incidence}, SIAM Journal on Imaging Sciences 8~(1) (2015)
  280--298.
\newline\urlprefix\url{https://doi.org/10.1137/140982787}

\bibitem{bao2015inverse}
G.~Bao, P.~Li, J.~Lin, F.~Triki,
  \href{https://doi.org/10.1088/0266-5611/31/9/093001}{Inverse scattering
  problems with multi-frequencies}, Inverse Problems 31~(9) (2015) 093001.
\newline\urlprefix\url{https://doi.org/10.1088/0266-5611/31/9/093001}

\bibitem{hanke1995landweber}
M.~Hanke, F.~Hettlich, O.~Scherzer,
  \href{https://doi.org/10.1115/DETC1995-0658}{The landweber iteration for an
  inverse scattering problem}, in: International Design Engineering Technical
  Conferences and Computers and Information in Engineering Conference, Vol.
  97669, American Society of Mechanical Engineers, 1995, pp. 909--915.
\newline\urlprefix\url{https://doi.org/10.1115/DETC1995-0658}

\bibitem{hettlich1999second}
F.~Hettlich, W.~Rundell, \href{https://doi.org/10.1137/S0036142998341246}{A
  second degree method for nonlinear inverse problems}, SIAM Journal on
  Numerical Analysis 37~(2) (1999) 587--620.
\newline\urlprefix\url{https://doi.org/10.1137/S0036142998341246}

\bibitem{hohage1997logarithmic}
T.~Hohage, \href{https://doi.org/10.1088/0266-5611/13/5/012}{Logarithmic
  convergence rates of the iteratively regularized gauss-newton method for an
  inverse potential and an inverse scattering problem}, Inverse problems 13~(5)
  (1997) 1279.
\newline\urlprefix\url{https://doi.org/10.1088/0266-5611/13/5/012}

\bibitem{kirsch1993domain}
A.~Kirsch, \href{https://doi.org/10.1088/0266-5611/9/1/005}{The domain
  derivative and two applications in inverse scattering theory}, Inverse
  problems 9~(1) (1993) 81.
\newline\urlprefix\url{https://doi.org/10.1088/0266-5611/9/1/005}

\bibitem{kress2003newton}
R.~Kress, \href{https://doi.org/10.1088/0266-5611/19/6/056}{Newton’s method
  for inverse obstacle scattering meets the method of least squares}, Inverse
  Problems 19~(6) (2003) S91.
\newline\urlprefix\url{https://doi.org/10.1088/0266-5611/19/6/056}

\bibitem{kress1994quasi}
R.~Kress, W.~Rundell, \href{https://doi.org/10.1088/0266-5611/10/5/011}{A
  quasi-newton method in inverse obstacle scattering}, Inverse Problems 10~(5)
  (1994) 1145.
\newline\urlprefix\url{https://doi.org/10.1088/0266-5611/10/5/011}

\bibitem{kress1997integral}
R.~Kress, Integral equation methods in inverse acoustic and electromagnetic
  scattering, Boundary integral formulations for inverse analysis(A 99-14241
  02-31), Southampton, United Kingdom and Billerica, MA, Computational
  Mechanics Publications (Advances in Boundary Elements Series) 1 (1997)
  67--92.

\bibitem{colton1996simple}
D.~Colton, A.~Kirsch, \href{https://doi.org/10.1088/0266-5611/12/4/003}{A
  simple method for solving inverse scattering problems in the resonance
  region}, Inverse problems 12~(4) (1996) 383.
\newline\urlprefix\url{https://doi.org/10.1088/0266-5611/12/4/003}

\bibitem{audibert2014generalized}
L.~Audibert, H.~Haddar, \href{https://doi.org/10.1088/0266-5611/30/3/035011}{A
  generalized formulation of the linear sampling method with exact
  characterization of targets in terms of farfield measurements}, Inverse
  Problems 30~(3) (2014) 035011.
\newline\urlprefix\url{https://doi.org/10.1088/0266-5611/30/3/035011}

\bibitem{kirsch1998characterization}
A.~Kirsch, \href{https://doi.org/10.1088/0266-5611/14/6/009}{Characterization
  of the shape of a scattering obstacle using the spectral data of the far
  field operator}, Inverse problems 14~(6) (1998) 1489.
\newline\urlprefix\url{https://doi.org/10.1088/0266-5611/14/6/009}

\bibitem{potthast2000stability}
R.~Potthast, \href{https://doi.org/10.1016/S0377-0427(99)00201-0}{Stability
  estimates and reconstructions in inverse acoustic scattering using singular
  sources}, Journal of computational and applied mathematics 114~(2) (2000)
  247--274.
\newline\urlprefix\url{https://doi.org/10.1016/S0377-0427(99)00201-0}

\bibitem{potthast2001point}
R.~Potthast, \href{https://doi.org/10.1201/9781420035483}{Point sources and
  multipoles in inverse scattering theory}, Chapman and Hall/CRC, 2001.
\newline\urlprefix\url{https://doi.org/10.1201/9781420035483}

\bibitem{ikehata1998reconstruction}
M.~Ikehata, \href{https://doi.org/10.1088/0266-5611/14/4/012}{Reconstruction of
  an obstacle from the scattering amplitude at a fixed frequency}, Inverse
  Problems 14~(4) (1998) 949.
\newline\urlprefix\url{https://doi.org/10.1088/0266-5611/14/4/012}

\bibitem{ikehata1999reconstruction}
M.~Ikehata, \href{https://doi.org/10.1016/S0165-2125(99)00006-2}{Reconstruction
  of obstacle from boundary measurements}, Wave motion 30~(3) (1999) 205--223.
\newline\urlprefix\url{https://doi.org/10.1016/S0165-2125(99)00006-2}

\bibitem{colton1998inverse}
D.~Colton, R.~Kress, \href{https://doi.org/10.1007/978-3-030-30351-8}{Inverse
  acoustic and electromagnetic scattering theory}, $4^{th}$ Edition, Vol.~93,
  Springer, 2019.
\newline\urlprefix\url{https://doi.org/10.1007/978-3-030-30351-8}

\bibitem{khoo2019switchnet}
Y.~Khoo, L.~Ying, \href{https://doi.org/10.1137/18M1222399}{Switchnet: a neural
  network model for forward and inverse scattering problems}, SIAM Journal on
  Scientific Computing 41~(5) (2019) A3182--A3201.
\newline\urlprefix\url{https://doi.org/10.1137/18M1222399}

\bibitem{yuwei2022solving}
Y.~Fan, L.~Ying, \href{https://dx.doi.org/10.4310/AMSA.2022.v7.n1.a2}{Solving
  inverse wave scattering with deep learning}, Annals of Mathematical Sciences
  and Applications 7~(1) (2022) 23--48.
\newline\urlprefix\url{https://dx.doi.org/10.4310/AMSA.2022.v7.n1.a2}

\bibitem{rekanos2002neural}
I.~T. Rekanos, \href{https://10.1109/20.996272}{Neural-network-based
  inverse-scattering technique for online microwave medical imaging}, IEEE
  transactions on magnetics 38~(2) (2002) 1061--1064.
\newline\urlprefix\url{https://10.1109/20.996272}

\bibitem{wei2018deep}
Z.~Wei, X.~Chen, \href{http://10.1109/TGRS.2018.2869221}{Deep-learning schemes
  for full-wave nonlinear inverse scattering problems}, IEEE Transactions on
  Geoscience and Remote Sensing 57~(4) (2018) 1849--1860.
\newline\urlprefix\url{http://10.1109/TGRS.2018.2869221}

\bibitem{adler2017solving}
J.~Adler, O.~{\"O}ktem, \href{https://10.1088/1361-6420/aa9581}{Solving
  ill-posed inverse problems using iterative deep neural networks}, Inverse
  Problems 33~(12) (2017) 124007.
\newline\urlprefix\url{https://10.1088/1361-6420/aa9581}

\bibitem{chen2019learning}
G.~Chen, P.~Shah, J.~Stang, M.~Moghaddam,
  \href{https://10.1109/TAP.2019.2948565}{Learning-assisted multimodality
  dielectric imaging}, IEEE Transactions on Antennas and Propagation 68~(3)
  (2019) 2356--2369.
\newline\urlprefix\url{https://10.1109/TAP.2019.2948565}

\bibitem{guo2019supervised}
R.~Guo, X.~Song, M.~Li, F.~Yang, S.~Xu, A.~Abubakar,
  \href{https://10.1109/TAP.2019.2902667}{Supervised descent learning technique
  for 2-d microwave imaging}, IEEE Transactions on Antennas and Propagation
  67~(5) (2019) 3550--3554.
\newline\urlprefix\url{https://10.1109/TAP.2019.2902667}

\bibitem{sanghvi2019embedding}
Y.~Sanghvi, Y.~Kalepu, U.~K. Khankhoje,
  \href{https://10.1109/TCI.2019.2915580}{Embedding deep learning in inverse
  scattering problems}, IEEE Transactions on Computational Imaging 6 (2019)
  46--56.
\newline\urlprefix\url{https://10.1109/TCI.2019.2915580}

\bibitem{chen2020review}
X.~Chen, Z.~Wei, M.~Li, P.~Rocca, \href{https://doi.org/10.2528/PIER20030705}{A
  review of deep learning approaches for inverse scattering problems (invited
  review)}, Progress In Electromagnetics Research 167 (2020) 67--81.
\newline\urlprefix\url{https://doi.org/10.2528/PIER20030705}

\bibitem{bao2007inverse}
G.~Bao, S.~Hou, P.~Li, \href{https://doi.org/10.1016/j.jcp.2007.08.020}{Inverse
  scattering by a continuation method with initial guesses from a direct
  imaging algorithm}, Journal of Computational Physics 227~(1) (2007) 755--762.
\newline\urlprefix\url{https://doi.org/10.1016/j.jcp.2007.08.020}

\bibitem{potthast1994frechet}
R.~Potthast, \href{https://dx.doi.org/10.1088/0266-5611/10/2/016}{Fr{\'e}chet
  differentiability of boundary integral operators in inverse acoustic
  scattering}, Inverse Problems 10~(2) (1994) 431.
\newblock \href {https://doi.org/10.1088/0266-5611/10/2/016}
  {\path{doi:10.1088/0266-5611/10/2/016}}.
\newline\urlprefix\url{https://dx.doi.org/10.1088/0266-5611/10/2/016}

\bibitem{kress1989linear}
R.~Kress, \href{https://doi.org/10.1007/978-1-4612-0559-3}{Linear integral
  equations}, Vol.~82, Springer, 1989.
\newline\urlprefix\url{https://doi.org/10.1007/978-1-4612-0559-3}

\bibitem{chandrasekaran2006fast}
S.~Chandrasekaran, P.~Dewilde, M.~Gu, W.~Lyons, T.~Pals,
  \href{https://doi.org/10.1137/050639028}{A fast solver for {HSS}
  representations via sparse matrices}, SIAM Journal on Matrix Analysis and
  Applications 29~(1) (2006) 67--81.
\newline\urlprefix\url{https://doi.org/10.1137/050639028}

\bibitem{chandrasekaran2006fast1}
S.~Chandrasekaran, M.~Gu, T.~Pals,
  \href{https://doi.org/10.1137/S0895479803436652}{A fast {ULV} decomposition
  solver for hierarchically semiseparable representations}, SIAM Journal on
  Matrix Analysis and Applications 28~(3) (2006) 603--622.
\newline\urlprefix\url{https://doi.org/10.1137/S0895479803436652}

\bibitem{gillman2012direct}
A.~Gillman, P.~M. Young, P.~G. Martinsson,
  \href{https://doi.org/10.1007/s11464-012-0188-3}{A direct solver with
  $\mathcal{O}({N})$ complexity for integral equations on one-dimensional
  domains}, Frontiers of Mathematics in China 7~(2) (2012) 217--247.
\newline\urlprefix\url{https://doi.org/10.1007/s11464-012-0188-3}

\bibitem{greengard2009fast}
L.~Greengard, D.~Gueyffier, P.~G. Martinsson, V.~Rokhlin,
  \href{https://doi.org/10.1017/S0962492906410011}{Fast direct solvers for
  integral equations in complex three-dimensional domains}, Acta Numerica 18
  (2009) 243--275.
\newline\urlprefix\url{https://doi.org/10.1017/S0962492906410011}

\bibitem{ho2012fast}
K.~L. Ho, L.~Greengard, \href{https://doi.org/10.1137/120866683}{A fast direct
  solver for structured linear systems by recursive skeletonization}, SIAM
  Journal on Scientific Computing 34~(5) (2012) A2507--A2532.
\newline\urlprefix\url{https://doi.org/10.1137/120866683}

\bibitem{martinsson2005fast}
P.~G. Martinsson, V.~Rokhlin,
  \href{https://doi.org/10.1016/j.jcp.2004.10.033}{A fast direct solver for
  boundary integral equations in two dimensions}, Journal of Computational
  Physics 205~(1) (2005) 1--23.
\newline\urlprefix\url{https://doi.org/10.1016/j.jcp.2004.10.033}

\bibitem{bebendorf2005hierarchical}
M.~Bebendorf, \href{https://doi.org/10.1007/s00607-004-0099-6}{Hierarchical
  {LU} decomposition-based preconditioners for {BEM}}, Computing 74~(3) (2005)
  225--247.
\newline\urlprefix\url{https://doi.org/10.1007/s00607-004-0099-6}

\bibitem{bormhierarchical}
S.~B{\"o}rm, L.~Grasedyck, W.~Hackbusch, Hierarchical matrices, Lecture notes
  21 (2003) 2003.

\bibitem{borm2003introduction}
S.~B{\"o}rm, L.~Grasedyck, W.~Hackbusch,
  \href{https://doi.org/10.1016/S0955-7997(02)00152-2}{Introduction to
  hierarchical matrices with applications}, Engineering Analysis with Boundary
  Elements 27~(5) (2003) 405--422.
\newline\urlprefix\url{https://doi.org/10.1016/S0955-7997(02)00152-2}

\bibitem{paige1982lsqr}
C.~C. Paige, M.~A. Saunders,
  \href{https://doi.org/10.1145/355984.355989}{{LSQR}: An algorithm for sparse
  linear equations and sparse least squares}, ACM Transactions on Mathematical
  Software (TOMS) 8~(1) (1982) 43--71.
\newline\urlprefix\url{https://doi.org/10.1145/355984.355989}

\bibitem{cakoni2005qualitative}
F.~Cakoni, D.~Colton, \href{https://doi.org/10.1007/3-540-31230-7}{Qualitative
  methods in inverse scattering theory: An introduction}, Springer Science \&
  Business Media, 2005.
\newline\urlprefix\url{https://doi.org/10.1007/3-540-31230-7}

\bibitem{barrow1977parametric}
H.~G. Barrow, J.~M. Tenenbaum, R.~C. Bolles, H.~C. Wolf, Parametric
  correspondence and chamfer matching: Two new techniques for image matching,
  in: Proceedings of the 5th International Joint Conference on Artificial
  Intelligence - Volume 2, IJCAI'77, Morgan Kaufmann Publishers Inc., San
  Francisco, CA, USA, 1977, p. 659–663.

\bibitem{Tygert2016}
M.~Tygert, J.~Bruna, S.~Chintala, Y.~LeCun, S.~Piantino, A.~Szlam,
  \href{https://doi.org/10.1162/neco_a_00824}{A mathematical motivation for
  complex-valued convolutional networks}, Neural Computation 28~(5) (2016)
  815--825.
\newblock \href {https://doi.org/10.1162/neco_a_00824}
  {\path{doi:10.1162/neco_a_00824}}.
\newline\urlprefix\url{https://doi.org/10.1162/neco_a_00824}

\bibitem{trabelsi2018deep}
C.~Trabelsi, O.~Bilaniuk, Y.~Zhang, D.~Serdyuk, S.~Subramanian, J.~F. Santos,
  S.~Mehri, N.~Rostamzadeh, Y.~Bengio, C.~J. Pal,
  \href{https://openreview.net/forum?id=H1T2hmZAb}{Deep complex networks}, in:
  International Conference on Learning Representations, 2018.
\newline\urlprefix\url{https://openreview.net/forum?id=H1T2hmZAb}

\bibitem{NEURIPS2019_9015}
A.~Paszke, S.~Gross, F.~Massa, A.~Lerer, J.~Bradbury, G.~Chanan, T.~Killeen,
  Z.~Lin, N.~Gimelshein, L.~Antiga, A.~Desmaison, A.~Kopf, E.~Yang, Z.~DeVito,
  M.~Raison, A.~Tejani, S.~Chilamkurthy, B.~Steiner, L.~Fang, J.~Bai,
  S.~Chintala,
  \href{https://proceedings.neurips.cc/paper_files/paper/2019/file/bdbca288fee7f92f2bfa9f7012727740-Paper.pdf}{Pytorch:
  An imperative style, high-performance deep learning library}, in: Advances in
  Neural Information Processing Systems 32, Curran Associates, Inc., 2019, pp.
  8024--8035.
\newline\urlprefix\url{https://proceedings.neurips.cc/paper_files/paper/2019/file/bdbca288fee7f92f2bfa9f7012727740-Paper.pdf}

\end{thebibliography}






\end{document}